\documentclass{amsart}

\newtheorem{theorem}{Theorem}[section]

\theoremstyle{definition}

\theoremstyle{remark}
\newtheorem{remark}[theorem]{Remark}

\numberwithin{equation}{section}

\usepackage{lipsum}
\usepackage{amsfonts}
\usepackage{graphicx}
\usepackage{epstopdf}

\usepackage{algorithm}
\usepackage{algorithmic}

\usepackage{amsopn}
\usepackage{bm}
\usepackage{multicol}
\usepackage{array}
\usepackage{tabularx}

\usepackage{graphicx}
\usepackage{blindtext}
\usepackage{changepage}
\usepackage{color}

\usepackage{caption,subcaption}
\usepackage{babel}
\usepackage{comment}
\usepackage{setspace}

\usepackage{array}
\newcolumntype{C}[1]{>{\centering\let\newline\\\arraybackslash\hspace{0pt}}m{#1}}

\usepackage{algorithmic}
\ifpdf
  \DeclareGraphicsExtensions{.eps,.pdf,.png,.jpg}
\else
  \DeclareGraphicsExtensions{.eps}
\fi

\newcommand{\RR}{\mathbb{R}}
\newcommand{\B}{\mathcal{B}}

\newcommand{\bk}{\mathbf{k}}
\newcommand{\bq}{\mathbf{q}}

\newcommand{\bj}{\mathbf{j}}

\newcommand{\bsigma}{\bm{\sigma}}

\newcommand{\D}{\mathcal{D}}

\newcommand{\n}{\bm{\mathrm{n}}}
\newcommand{\bl}{\bm{\mathrm{l}}}
\newcommand{\m}{\bm{\mathrm{m}}}

\newcommand*{\im}{\mathop{}\!\mathrm{i}}
\newcommand*{\e}{\mathop{}\!\mathrm{e}}

\newcommand{\rd}{\mathrm{d}}

\begin{document}

\title{A fast Fourier spectral method for wave kinetic equation}

\author{Kunlun Qi}
\address{Department of Computational Mathematics, Science and Engineering and Department of Mathematics, Michigan State University, East Lansing, MI 48824}
\email{qikunlun@msu.edu, kunlunqi.math@gmail.com}

\author{Lian Shen}
\address{Department of Mechanical Engineering, University of Minnesota--Twin Cities, Minneapolis, MN 55455}
\email{shen@umn.edu}
\thanks{L.S. is funded in part by ONR and the University of Minnesota.}

\author{Li Wang}
\address{School of Mathematics, University of Minnesota--Twin Cities, Minneapolis, MN 55455}
\email{liwang@umn.edu}
\thanks{L.W. is funded in part by NSF grant DMS-1846854 and the Simons Fellowship.}


\date{Submitted to the editors DATE}

\emergencystretch 3em


\begin{abstract}
The central object in wave turbulence theory is the wave kinetic equation (WKE), which is an evolution equation for wave action density and acts as the wave analog of the Boltzmann kinetic equations for particle interactions. Despite recent exciting progress in the theoretical aspects of the WKE, numerical developments have lagged behind. In this paper, we introduce a fast Fourier spectral method for solving the WKE. The key idea lies in reformulating the high-dimensional nonlinear wave kinetic operator as a spherical integral, analogous to the classical Boltzmann collision operator. The conservation of mass and momentum leads to a double convolution structure in Fourier space, which can be efficiently handled by the fast Fourier transform (FFT), reducing the computational cost from $O(N^{3d})$ to $O(M N^d \log N)$ with $N$-frequency nodes and $M \ll N^{2d-1}$ in $d$ dimensions. We demonstrate the accuracy and efficiency of the proposed method through several numerical tests in both 2D and 3D, revealing and conjecturing some interesting and unique features of this equation.
\end{abstract}

\maketitle

\section{Introduction}
\label{sec:intro}


In the past decades, \textit{wave turbulence theory} (WTT) has been developed as a framework for understanding the statistical behavior of nonlinear wave systems in non-equilibrium settings \cite{zakharov2012kolmogorov, nazarenko2011wave, newell2011wave}, encompassing phenomena such as ocean surface gravity waves \cite{janssen2004interaction}, capillary water waves \cite{zakharov1967weak}, nonlinear optics \cite{dyachenko1992optical} and so forth. As its fundamental governing equation, \textit{wave kinetic equation} (WKE), is derived as the kinetic limit of time-reversible wave dispersive equations with weakly nonlinear interaction (i.e., waves of different frequencies interact nonlinearly at the microscopic scale) \cite{CG2019, CG2020, DH2021, DH2023}, describing the evolution of energy densities.

There have been notable advances in the theoretical understanding of WKE. 
The well-posedness of WKE, for the  Schr{\"o}dinger dispersion relation $\omega(\bk) = |\bk|^2$ as in \eqref{WKE}, was proved by Escobedo and Vel\'azquez in their pioneering series of work \cite{EV2015-1, EV2015-2}, where they proved the existence of global measure-valued solutions in the isotropic case. Furthermore, Germain, Ionescu, and Tran extended the local well-posedness for general dispersion relations in \cite{GIT2020}.
Additionally, the energy cascade and condensation of the solution were demonstrated by Staffilani and Tran in \cite{ST2024, ST2024-2}.
The rigorous derivation of the WKE has seen rapid advancement since the seminal work of Lukkarinen and Spohn \cite{LS2011}. In particular, Deng and Hani, in their recent series of papers \cite{DH2021, DH2023, DH-long}, provided a rigorous derivation of the homogeneous WKE from the cubic nonlinear Schr{\"o}dinger equation over various kinetic timescales. Additionally, we refer to the work of Buckmaster, Germain, Hani, and Shatah \cite{BGHS2020, BGHS2021}, and references therein, for further progress on the rigorous derivation of the WKE.

Although the WKE has demonstrated broad applicability in real-world problems and has seen substantial theoretical progress, its numerical simulation has lagged behind. This delay is due to the challenges posed by the equation’s complex structure, which is characterized by high dimensionality and nonlinearities.
The numerical study of the WKE dates back to the 1980s, with two dominant approaches being the Discrete Interaction Approximation (DIA) \cite{hasselmann1985computations} and quadrature-based methods (known as the WRT method) \cite{webb1978non, tracy1982theory, resio1991numerical}. These methods mentioned above either have limited validity, capturing only specific aspects of the phenomena, or are too computationally expensive to be practical. Little progress was made until recently, when Walton and Tran applied a finite volume scheme to simulate the 3-WKE \cite{WT2023}, a variant of the WKE with fewer wave interactions. Additionally, deep learning techniques have been employed to address the numerical challenges of solving the WKE in \cite{WTB2024}. Other related developments can be found in \cite{BBKKS2022, SMNF2024} and the references therein.


On the other hand, as one of the most efficient and accurate deterministic type methods, the Fourier spectral method has shown its success and provided an effective framework in approximating the classical kinetic equation, i.e., Boltzmann equation, as first established in \cite{PP96, PR00}. This method offers several advantages: (i) it is deterministic and yields highly accurate results compared to stochastic methods; (ii) it exploits the translation-invariant nature of the Boltzmann collision operator using Fourier bases; and (iii) upon Galerkin projection, the collision operator adopts a convolution-like structure, enabling further acceleration via the fast Fourier transform (FFT) \cite{MP06, GHHH17, HQ2020}. Due to these merits, the Fourier spectral method has gained significant popularity over the past decade for solving the Boltzmann equation and related collisional kinetic models. For examples of their application, see \cite{PRT00, FR03, FMP06, HY12, FHJ12, HQY20}, as well as the review in \cite{pareschi}.


\textbf{Motivation and our contribution:} Given the discussions above, this paper extends the Fourier spectral method to solve the WKE in its complete form, for both isotropic and non-isotropic cases. Additionally, we propose a fast algorithm to improve computational efficiency. The key idea behind applying the Fourier spectral method to the wave kinetic operator \eqref{Kre} is to reformulate the high-dimensional nonlinear operator as an integral over the sphere, similar to the classical Boltzmann collision operator.
To further accelerate the computation, we exploit the double-convolution structure inherent in the weighted summation numerical system through specific quadrature and reorganization. This structure allows the use of the Fast Fourier Transform (FFT), which reduces the computational cost from $O(N^{3d})$ to $O(M N^d \log N)$, where $N$ represents the number of frequency nodes and $M \ll N^{2d-1}$ in $d$ dimensions. Notably, all evaluations are performed ``on the fly", eliminating the need for pre-computation or extensive storage, further enhancing the practicality and scalability of the algorithm.


\textbf{Organization of the paper:} The rest of this paper is organized as follows. The WKE and its associated physical quantities and stationary states are introduced in Section \ref{sec:wke}. In Section \ref{sec:spectral}, we propose the numerical formulation when applying the Fourier spectral method to the WKE, as well as the fast version by FFT. Numerical examples are presented in 2D and 3D to validate our proposed method in Section \ref{sec:numerical}. Some conclusion remarks are given in Section \ref{sec:conclusion}.

\section{Wave kinetic equation}
\label{sec:wke}

The 4-wave kinetic equation (4-WKE) reads 
\begin{equation}\label{WKE}
    \partial_t f(t,\bk) = \mathcal{K}(f,f,f)(t,\bk), \quad (t,\bk) \in \RR^{+} \times \RR^d,
\end{equation}
where
\begin{equation}\label{Kre}
\begin{split}
    \mathcal{K}(f,f,f)&(\bk) = \int_{(\RR^d)^3} \frac{f(\bk) f(\bk_1) f(\bk_2) f(\bk_3)}{|\bk|^{\frac{\beta}{2}} |\bk_1|^{\frac{\beta}{2}} |\bk_2|^{\frac{\beta}{2}} |\bk_3|^{\frac{\beta}{2}}} \left[\frac{1}{f(\bk)} - \frac{1}{f(\bk_1)} + \frac{1}{f(\bk_2)} - \frac{1}{f(\bk_3)}\right]\\[4pt]
    &\times \delta(\bk_1- \bk_2 + \bk_3 - \bk) \delta \Big( \omega(\bk_1) - \omega(\bk_2) + \omega(\bk_3) - \omega(\bk) \Big) \,\rd \bk_1 \,\rd \bk_2 \,\rd \bk_3\,.
\end{split}
\end{equation}
Here, $\bk \mapsto \omega(\bk)$ is the dispersion relation, and the parameter $\beta$ in the cross-section depends on the nonlinearity of the dynamic system at the microscopic scale. In particular, $\beta = 0$ corresponds to the cubically nonlinear Schr{\"o}dinger equation \cite{EV2015-1, EV2015-2, Menegaki2024}. $\delta$ here is the Dirac Delta function. 


Although radial symmetry is typically assumed in theoretical analysis, this paper considers the general function $f(t,\bk)$, where $ \bk \in \mathbb{R}^d, d=2,3$, and focuses particularly on the dispersion relation 
$\omega(\bk) = \frac{|\bk|^2}{2}$ with $\beta =0$. This particular case allows for a rigorous derivation of the wave kinetic equation (WKE) from the Schrodinger equation \cite{DH2023}. The same methodology developed here can be applied to water gravity waves, where the dispersion relation is $\omega(\bk) = \sqrt{g|\bk|}$ with $g$ being the gravitational acceleration, with only minor modifications. Applications to other dispersion relations, such as those for ocean surface gravity, capillary, and gravity-capillary waves, will be explored in future work.

Then $\mathcal{K}(f,f,f)$ in \eqref{Kre} takes the following form: 
\begin{equation}\label{K}
\begin{split}
    \mathcal{K}(f,f,f)(\bk)& = \int_{(\RR^d)^3} \left[ f_1 f_2 f_3 - f f_2 f_3 + f_1 f f_3 - f_1 f_2 f \right]\\
    &\times \delta(\bk_1- \bk_2 + \bk_3 - \bk) \delta \left( \frac{|\bk_1|^2}{2} - \frac{|\bk_2|^2}{2} + \frac{|\bk_3|^2}{2} - \frac{|\bk|^2}{2}\right) \,\rd \bk_1 \,\rd \bk_2 \,\rd \bk_3.
\end{split}
\end{equation}
where $f:=f(t,\bk)$ and $f_i := f(t,\bk_i), i=1,2,3$.

Let 
\begin{equation}
    \rho := \int_{\RR^d} f(\bk)\,\rd \bk, \quad  \quad E := \int_{\RR^d} f(\bk) |\bk|^2 \,\rd \bk,
\end{equation}
be the macroscopic density and energy, respectively. Then, their conservation can be  derived from the weak formulation of \eqref{K}:
\begin{equation}\label{weak-WKE}
\begin{split}
    \int_{\RR^d} \mathcal{K}(f,f,f)(\bk) \varphi(\bk) \,\rd \bk& = \frac{1}{4}\int_{(\RR^d)^4}  f f_1 f_2 f_3  \left[ \frac{1}{f} - \frac{1}{f_1} + \frac{1}{f_2} - \frac{1}{f_3}\right]\\[4pt]
    &\times \delta(\bk_1- \bk_2 + \bk_3 - \bk) \times \delta \left( \frac{|\bk_1|^2}{2} - \frac{|\bk_2|^2}{2} + \frac{|\bk_3|^2}{2} - \frac{|\bk|^2}{2}\right) \\
    & \times \Big[ \varphi(\bk) +  \varphi(\bk_2) - \varphi(\bk_1) - \varphi(\bk_3)\Big] \,\rd \bk_1 \,\rd \bk_2 \,\rd \bk_3 \,\rd \bk\,,
\end{split}
\end{equation}
where $\varphi$ is a test function. Indeed, the conservation of mass and energy corresponds to choosing $\varphi$ to be $1$ and $|\bk|^2$, respectively. 

Moreover, from \eqref{Kre}, one observes that when setting 
\begin{equation}\label{equi}
    f_{\text{eq}}:= \frac{1}{\mu + \bm{\nu} \cdot \bk + \xi |\bk|^2}\,,
\end{equation}
where $(\mu, \bm{\nu}, \xi) \in \RR \times \RR^d \times \RR$ are such that $\mu + \bm{\nu} \cdot \bk + \xi |\bk|^2 > 0 $ for any $\bk \in \RR^d$ ($d=2, 3$), 
we have 
$\mathcal K (f_{\text{eq}}, f_{\text{eq}}, f_{\text{eq}}) \equiv 0$. This is the so-called Rayleigh-Jeans (RJ) stationary solution \cite{GIT2020, Menegaki2024}. This solution is shown to have local $L^2$-stability \cite{Menegaki2024}, with the only instability occurs when $\mu = 0, \bm{\nu} =\bm{0}$ \cite{EM2024}.

Additionally, \eqref{WKE} admits two Kolmogorov-Zakharov (KZ) steady states, which are associated with the cascade phenomena. Specifically, in the isotropic case where $f$ 
depends only on $|\bk|$,  $f_{\text{eq}} = \omega^{-\frac{7}{6}}(\bk)$ corresponds to an inverse cascade of mass, where mass is transferred from infinite to zero frequencies; while $f_{\text{eq}} = \omega^{-\frac{3}{2}}(\bk)$ corresponds to a direct energy cascade, where energy is transferred from zero to infinite frequencies \cite{BZ1998, EV2015-1, EV2015-2}. For further discussions on the stability of the KZ solutions, we refer to the recent work by Collot, Dietert, and Germain in \cite{CDG2024}.

To facilitate the numerical computation, we reformulate the formula \eqref{K} into an integral over the sphere centered at $\frac{\bk+\bk_2}{2}$ with radius $\frac{|\bk-\bk_2|}{2}$ by introducing $\bsigma \in \mathbb{S}^{d-1}$, as in the classical Boltzmann-type collision operator \cite{HY15} (details are provided in Appendix \ref{appen:derivation}):
\begin{equation}\label{WKE-Boltz}
    \begin{split}
    \mathcal{K}(f,f,f) =& \frac{1}{2^{d-1}} \int_{\RR^d} \int_{\mathbb{S}^{d-1}}  \Big[ f\left( \bk' \right) f(\bk_2) f\left(\bk'_2\right)  - f(\bk) f(\bk_2) f\left(\bk'_2\right) + f\left(\bk'\right) f(\bk) f\left(\bk'_2\right)\\[6pt]
    &\qquad \qquad \qquad - f\left(\bk'\right) f(\bk_2) f(\bk)\Big] |\bk-\bk_2|^{d-2} \,\rd \bsigma \,\rd \bk_2\\[3pt]
    :=& \frac{1}{2^{d-1}} \Big[ \mathcal{K}_1(f,f,f) - \mathcal{K}_2(f,f,f) + \mathcal{K}_3(f,f,f) - \mathcal{K}_4(f,f,f)\Big]\,,
    \end{split}
\end{equation}
where
\begin{equation}\label{bkbk'}
    \bk' = \frac{\bk_2 + \bk}{2} + \frac{|\bk-\bk_2|}{2}\bsigma, \quad 
    \bk'_2 = \frac{\bk_2 + \bk}{2} - \frac{|\bk-\bk_2|}{2}\bsigma\,,
\end{equation}
and
\begin{equation*}
\left\{
\begin{aligned}
    \mathcal{K}_1(f,f,f)(\bk):=& \int_{\RR^d} \int_{\mathbb{S}^{d-1}} |\bk-\bk_2|^{d-2} \, f\left( \bk' \right) f(\bk_2) f\left(\bk'_2\right)  \,\rd \bsigma \,\rd \bk_2,\\[3pt]
    \mathcal{K}_2(f,f,f)(\bk):=&  \int_{\RR^d} \int_{\mathbb{S}^{d-1}}  |\bk-\bk_2|^{d-2} \, f(\bk) f(\bk_2) f\left(\bk'_2\right)  \,\rd \bsigma \,\rd \bk_2,\\[3pt]
    \mathcal{K}_3(f,f,f)(\bk):=&  \int_{\RR^d} \int_{\mathbb{S}^{d-1}} |\bk-\bk_2|^{d-2} \, f\left(\bk'\right) f(\bk) f\left(\bk'_2\right) \,\rd \bsigma \,\rd \bk_2,\\[3pt]
    \mathcal{K}_4(f,f,f)(\bk):=& \int_{\RR^d} \int_{\mathbb{S}^{d-1}} |\bk-\bk_2|^{d-2} \, f\left(\bk'\right) f(\bk_2) f(\bk) \,\rd \bsigma \,\rd \bk_2.
\end{aligned}
\right.
\end{equation*}
To close this section, we would like to highlight the similarities and differences between our approach and two previous methods that inspired our work. The first, developed in \cite{HY15}, addresses the energy-space boson Boltzmann equation. In this work, the collision kernel depends on the minimum of the four interacting energies, and a special domain decomposition design is proposed to handle this interaction kernel. However, their approach is limited to the isotropic case, whereas we extend the analysis to the general case.
The second method, developed in \cite{HY12}, applies to the quantum Boltzmann operator. Unlike their approach, which utilizes the Carleman representation following \cite{MP06}, our method does not rely on this representation. This makes our approach more flexible, as it can be extended to cases with different interaction kernels.

\section{Fourier spectral method for the 4-WKE}
\label{sec:spectral}

In this section, we present the complete numerical approximation procedure used to solve the 4-WKE. For clarity, we first provide a global overview of the
numerical pipeline, and then describe in detail each component of the Fourier spectral method, including truncation, discretization, kernel approximation, and fast evaluation of the collision operator. We conclude with an explicit analysis of the computational complexity,
highlighting the efficiency of the proposed algorithm.

\subsection{Overview of the numerical procedure}
\label{subsec:overview-numerics}

The numerical method developed in this work consists of the following main steps:
\begin{enumerate}
  \item Truncation of the wave-number domain $\bk$ and periodic extension of the solution;
  \item Fourier-Galerkin discretization in the wave-number variable $\bk$;
  \item Reformulation of the 4-wave collision operator $\mathcal{K}(f,f,f)$ into a Boltzmann-type integral revealing a convolution structure as in \eqref{WKE-Boltz};
  \item Approximation of the pre-computed weight functions that consists of the collision kernel by a low-rank separated representation using radial and angular quadratures;
  \item Fast evaluation of the resulting double-convolutional structure using FFTs.
\end{enumerate}
This organization is intended to clearly separate the standard discretization effects from fast algorithm effects, and to make transparent the origin of the computational efficiency of the method.

\subsection{Numerical formulation}
\label{subsec:formulation}

Now, we first present a step-by-step description of the numerical formulation and discretization process.

\noindent\textit{Step 1: Domain truncation and periodization.}\\
To apply the Fourier spectral method, we assume that $f(\bk)$ is compactly supported in $\mathcal{B}_S$, i.e., $\text{supp} \big(f(\bk)\big) \subset \mathcal{B}_S$, where $\mathcal{B}_S$ is the ball of radius $S$ centered in the origin. By conservation of energy, $\bk'$ defined in \eqref{bkbk'} satisfies: 
\[
(\bk')^2 \leq \bk^2 + \bk_2^2 \leq 2S^2,
\]
and similarly, $(\bk'_2)^2 \leq 2S^2$. Therefore, $\text{supp}\big( \mathcal{K}(f,f,f)\big) \subset \mathcal{B}_{\sqrt{2}S}$.
Next, we consider the change of variable $\bk \mapsto \bq := \bk-\bk_2 \in \mathcal{B}_R$, where $R = 2S$ since $|\bq| = |\bk-\bk_2|\leq 2S$. We restrict the function $f(\bq)$ to the larger domain $\D_L = [-L, L]^d$ and extend it periodically,  where $L$ is chosen to satisfy $L \geq \frac{3+\sqrt{2}}{2} S$ to avoid aliasing errors \cite{pareschi}.

We want to point out that this cut-off assumption is not physically motivated, but is made primarily for numerical purposes as for the classical Boltzmann equation \cite{pareschi}. Additionally, it aids in the study of the long-time behavior and the existence of solutions around the Rayleigh-Jeans equilibrium solutions \cite{Menegaki2024}.

\begin{remark}
    We note that unlike the exponentially decay Maxwellian equilibrium of the Boltzmann equation \cite{HQ2020}, the RJ equilibrium \eqref{equi} of 4-WKE exhibits algebraic decay. Nevertheless, domain truncation remains justified in our context to certain degree. For instance, in the isotropic case, the dynamics are dominated by the Dirac measure type blow-up at the origin \cite[Page 2-6]{EV2015-1} (analogous to Bose-Einstein condensation), resulting in negligible tails. For cases involving heavier tails, such as perturbations near the RJ equilibrium \cite{Menegaki2024}, we mitigate aliasing by selecting a domain size $L$ sufficiently large to ensure boundary values remain small. Additionally, the use of classical de-aliasing methods (e.g., the Orszag rule \cite{Orszag1971-1, PO1971} or Hou-Li filter \cite{HL2007}), facilitated by our fast algorithm, can also ensure that the truncation error remains controlled.
\end{remark}

Then, $\mathcal{K}_1$ can be written as follows: 
\begin{equation*}
\begin{split}
    \mathcal{K}_1^R(f,f,f)(t,\bk) = \int_{\mathcal{B}_R} \int_{\mathbb{S}^{d-1}} |\bq|^{d-2} \, f\left( \bk' \right) f(\bk-\bq) f\left(\bk'_2\right)  \,\rd \bsigma \,\rd \bq
\end{split}
\end{equation*}
with
\begin{equation*}
    \bk' = \bk - \frac{1}{2}(\bq - |\bq|\bsigma), \quad 
    \bk'_2 = \bk - \frac{1}{2}(\bq + |\bq|\bsigma)\,.
\end{equation*}

\noindent\textit{Step 2: Spectral approximation.} \\
We approximate $f$ by a truncated Fourier series ($\bj$ is the $d$-dimensional Fourier mode):
\begin{equation}\label{fj}
    f(\bk) \approx f_N(\bk):= \sum_{|\bj| = -\frac{N}{2}}^{\frac{N}{2}-1} \hat{f}_{\bj} \e^{\im \frac{\pi}{L} \bj \cdot \bk}, \quad \text{with}\quad \hat{f}_{\bj} = \frac{1}{(2L)^{d}} \int_{\D_L} f(\bk) \e^{-\im \frac{\pi}{L} \bj \cdot \bk} \,\rd \bk.
\end{equation}

\noindent\textit{Step 3: Galerkin projection.} \\
Substituting $\eqref{fj}$ into \eqref{WKE} and projecting onto the space spanned by the orthogonal basis $\e^{\im \frac{\pi}{L}\bj \cdot \bk}, \, \bj \in \mathbb{Z}^d$, we obtain the following evolution equation for $\hat{f}_{\bj}$:
\begin{equation}\label{Kj}
    \frac{\rd \hat{f}_{\bj}}{\rd t} = \hat{\mathcal{K}}_{\bj}^R = \hat{\mathcal{K}}_{1,\bj}^R - \hat{\mathcal{K}}_{2,\bj}^R + \hat{\mathcal{K}}_{3,\bj}^R - \hat{\mathcal{K}}_{4,\bj}^R \,,
\end{equation}
where
\begin{equation*}
    \hat{\mathcal{K}}_{\ell,\bj}^R := \frac{1}{(2L)^{d}} \int_{\D_L} \mathcal{K}_\ell^R(\bk) \e^{-\im \frac{\pi}{L} \bj \cdot \bk} \,\rd \bk, \quad \text{for} \ \ell=1,2,3,4. 
\end{equation*}

Then, by substituting \eqref{fj} into $\hat{\mathcal{K}}_{1,\bj}$, we get
\begin{equation}\label{KK1}
    \begin{split}
        \hat{\mathcal{K}}_{1,\bj}^R
        = & \frac{1}{(2L)^{d}} \int_{\D_L} \int_{\B_R} \int_{\mathbb{S}^{d-1}} |\bq|^{d-2}  f_N(\bk') f_N(\bk-\bq) f_N(\bk'_2)  \,\rd \bsigma \,\rd \bq \e^{-\im \frac{\pi}{L} \bj \cdot \bk} \,\rd \bk\\[6pt]
        =& \frac{1}{(2L)^{d}} \int_{\D_L} \int_{\B_R} \int_{\mathbb{S}^{d-1}} |\bq|^{d-2}
        \left( \sum_{|\bl| = -\frac{N}{2}}^{\frac{N}{2}-1} \hat{f}_{\bl} \e^{\im \frac{\pi}{L} \bl \cdot \bk'} \right)
        \left(\sum_{|\m| = -\frac{N}{2}}^{\frac{N}{2}-1} \hat{f}_{\m} \e^{\im \frac{\pi}{L} \m \cdot (\bk-\bq)}\right)\\
        &\times \left(\sum_{|\n| = -\frac{N}{2}}^{\frac{N}{2}-1} \hat{f}_{\n} \e^{\im \frac{\pi}{L} \n \cdot \bk'_2}\right) 
         \e^{-\im \frac{\pi}{L} \bj \cdot \bk} \,\rd \bsigma \,\rd \bq  \,\rd \bk\\[6pt]
        =& \frac{1}{(2L)^{d}} \sum_{|\bl| = -\frac{N}{2}}^{\frac{N}{2}-1} \sum_{|\m| = -\frac{N}{2}}^{\frac{N}{2}-1} \sum_{|\n| = -\frac{N}{2}}^{\frac{N}{2}-1} \hat{f}_{\bl} \hat{f}_{\m} \hat{f}_{\n} \int_{\D_L} \int_{\B_R} \int_{\mathbb{S}^{d-1}} |\bq|^{d-2}\\
        &\times \left(  \e^{\im \frac{\pi}{L} \bl \cdot [\bk - \frac{1}{2}(\bq - |\bq|\bsigma)])} \right)
        \left(  \e^{\im \frac{\pi}{L} \m \cdot (\bk-\bq)}\right)
        \left(  \e^{\im \frac{\pi}{L} \n \cdot [\bk - \frac{1}{2}(\bq + |\bq|\bsigma)]}\right) 
        \e^{-\im \frac{\pi}{L} \bj \cdot \bk} \,\rd \bsigma \,\rd \bq  \,\rd \bk\\[6pt]
        =& \frac{1}{(2L)^{d}} \sum_{|\bl| = -\frac{N}{2}}^{\frac{N}{2}-1} \sum_{|\m| = -\frac{N}{2}}^{\frac{N}{2}-1} \sum_{|\n| = -\frac{N}{2}}^{\frac{N}{2}-1} \hat{f}_{\bl} \hat{f}_{\m} \hat{f}_{\n} \int_{\D_L} \int_{\B_R} \int_{\mathbb{S}^{d-1}} |\bq|^{d-2}\\
        &\times   \e^{\im \frac{\pi}{L} (\bl+\m+\n-\bj) \cdot \bk}
          \e^{-\im \frac{\pi}{L} \frac{1}{2} (2\m+\bl+\n) \cdot \bq} 
          \e^{\im \frac{\pi}{L} \frac{1}{2}|\bq|(\bl-\n)\cdot\bsigma}  
          \,\rd \bsigma \,\rd \bq \,\rd \bk\\[6pt]
        =& \sum_{|\bl| = -\frac{N}{2}}^{\frac{N}{2}-1} \sum_{|\m| = -\frac{N}{2}}^{\frac{N}{2}-1} \sum_{|\n| = -\frac{N}{2}}^{\frac{N}{2}-1} \hat{f}_{\bl} \hat{f}_{\m} \hat{f}_{\n}  \int_{\B_R} \int_{\mathbb{S}^{d-1}} |\bq|^{d-2}\\
        &\times   \delta(\bl+\m+\n-\bj) 
          \e^{-\im \frac{\pi}{L} \frac{1}{2} (2\m+\bl+\n) \cdot \bq} 
          \e^{\im \frac{\pi}{L} \frac{1}{2}|\bq| (\bl-\n)\cdot\bsigma} \,\rd \bsigma \,\rd \bq  \\[6pt]
        =& \sum_{\substack{|\bl|,|\m|, |\n|= -\frac{N}{2},\\ \bl+\n=\bj-\m}}^{\frac{N}{2}-1}  G_1(\bl,\m,\n)\hat{f}_{\bl} \hat{f}_{\m} \hat{f}_{\n}  \,.
    \end{split}
\end{equation}
Here $\delta$ is a Kronecker delta function, i.e., $\delta (n) = 1$ when $n = 0$ and zero otherwise , and the weight function $G_1(\bl,\m,\n)$ takes the form:
\begin{equation*}
    G_1(\bl,\m,\n) := \int_{\B_R} \int_{\mathbb{S}^{d-1}} |\bq|^{d-2}
        \e^{-\im \frac{\pi}{L} \frac{1}{2} (2\m+\bl+\n) \cdot \bq} 
        \e^{\im \frac{\pi}{L} \frac{1}{2}|\bq| (\bl-\n)\cdot\bsigma} \,\rd \bsigma \,\rd \bq \,.
\end{equation*}

At this point, we can summarize the computational cost of directly evaluating the operator $\mathcal{K}_1^R$ using the classical spectral method with the form \eqref{KK1}:
\begin{itemize}
    \item Pre-compute the weight function $G_1(\bl,\m,\n)$ -- storage requirement $O(N^{3d})$;
    \item Compute the $\hat{f}_{\bj}$ by using the FFT -- cost $O(N^d \log N)$;
    \item Compute the $\hat{\mathcal{K}}_{1,\bj}^R$ in the weighted summation form \eqref{KK1} -- cost $O(N^{3d})$;
    \item Compute the $\mathcal{K}_1^R$ by using the inverse FFT to $\hat{\mathcal{K}}_{1,\bj}^R$ -- cost $O(N^d \log N)$.
\end{itemize}

Note that $\hat{\mathcal{K}}_{1,\bj}^R$ is the most complicated one, the rest $\hat{\mathcal{K}}_{2,\bj}^R, \hat{\mathcal{K}}_{3,\bj}^R, \hat{\mathcal{K}}_{4,\bj}^R$ can be derived in a similar but simpler manner as follows. 

For $\hat{\mathcal{K}}_{2,\bj}$, by substituting \eqref{fj} into $\hat{\mathcal{K}}_{2,\bj}$, we have
\begin{equation}\label{KK2}
    \begin{split}
        \hat{\mathcal{K}}_{2,\bj}^R 
        =& \frac{1}{(2L)^{d}} \int_{\D_L} \int_{\B_R} \int_{\mathbb{S}^{d-1}} |\bq|^{d-2}
        \left( \sum_{|\bl| = -\frac{N}{2}}^{\frac{N}{2}-1} \hat{f}_{\bl} \e^{\im \frac{\pi}{L} \bl \cdot \bk} \right)
        \left(\sum_{|\m| = -\frac{N}{2}}^{\frac{N}{2}-1} \hat{f}_{\m} \e^{\im \frac{\pi}{L} \m \cdot (\bk-\bq)}\right)\\
        & \qquad \times \left(\sum_{|\n| = -\frac{N}{2}}^{\frac{N}{2}-1} \hat{f}_{\n} \e^{\im \frac{\pi}{L} \n \cdot \bk'_2}\right) 
         \e^{-\im \frac{\pi}{L} \bj \cdot \bk} \,\rd \bsigma \,\rd \bq  \,\rd \bk\\[6pt]
        =& \frac{1}{(2L)^{d}} \sum_{|\bl| = -\frac{N}{2}}^{\frac{N}{2}-1} \sum_{|\m| = -\frac{N}{2}}^{\frac{N}{2}-1} \sum_{|\n| = -\frac{N}{2}}^{\frac{N}{2}-1} \hat{f}_{\bl} \hat{f}_{\m} \hat{f}_{\n} \int_{\D_L} \int_{\B_R} \int_{\mathbb{S}^{d-1}} |\bq|^{d-2}\\
        &\times   \e^{\im \frac{\pi}{L} \bl \cdot \bk}
        \left(  \e^{\im \frac{\pi}{L} \m \cdot (\bk-\bq)}\right)
        \left(  \e^{\im \frac{\pi}{L} \n \cdot [\bk - \frac{1}{2}(\bq + |\bq|\bsigma)]}\right) 
        \e^{-\im \frac{\pi}{L} \bj \cdot \bk} \,\rd \bsigma \,\rd \bq  \,\rd \bk\\[4pt]
        =& \sum_{\substack{|\bl|,|\m|, |\n|= -\frac{N}{2},\\ \bl+\n=\bj-\m}}^{\frac{N}{2}-1}  G_2(\m,\n)\hat{f}_{\bl} \hat{f}_{\m} \hat{f}_{\n} \,,
    \end{split}
\end{equation}
where 
\begin{equation*}
    G_2(\m,\n) : = \int_{\B_R} \int_{\mathbb{S}^{d-1}} |\bq|^{d-2}
        \e^{-\im \frac{\pi}{L} \frac{1}{2} (2\m+\n) \cdot \bq} 
        \e^{-\im \frac{\pi}{L} \frac{1}{2}|\bq| \n\cdot\bsigma} \,\rd \bsigma \,\rd \bq.
\end{equation*}

For $\hat{\mathcal{K}}_{3,\bj}^R$, 
\begin{equation}\label{KK3}
    \begin{split}
        \hat{\mathcal{K}}_{3,\bj}^R
        =& \frac{1}{(2L)^{d}} \int_{\D_L} \int_{\B_R} \int_{\mathbb{S}^{d-1}} |\bq|^{d-2}
        \left( \sum_{|\bl| = -\frac{N}{2}}^{\frac{N}{2}-1} \hat{f}_{\bl} \e^{\im \frac{\pi}{L} \bl \cdot \bk'} \right)
        \left(\sum_{|\m| = -\frac{N}{2}}^{\frac{N}{2}-1} \hat{f}_{\m} \e^{\im \frac{\pi}{L} \m \cdot \bk}\right)\\
        &\times
        \left(\sum_{|\n| = -\frac{N}{2}}^{\frac{N}{2}-1} \hat{f}_{\n} \e^{\im \frac{\pi}{L} \n \cdot \bk'_2}\right)
         \e^{-\im \frac{\pi}{L} \bj \cdot \bk} \,\rd \bsigma \,\rd \bq  \,\rd \bk\\[6pt]
        =& \frac{1}{(2L)^{d}} \sum_{|\bl| = -\frac{N}{2}}^{\frac{N}{2}-1} \sum_{|\m| = -\frac{N}{2}}^{\frac{N}{2}-1} \sum_{|\n| = -\frac{N}{2}}^{\frac{N}{2}-1} \hat{f}_{\bl} \hat{f}_{\m} \hat{f}_{\n} \int_{\D_L} \int_{\B_R} \int_{\mathbb{S}^{d-1}} |\bq|^{d-2}\\
        &\times   \e^{\im \frac{\pi}{L} (\bl+\m+\n-\bj) \cdot \bk}
          \e^{-\im \frac{\pi}{L} \frac{1}{2} (\bl+\n) \cdot \bq} 
          \e^{\im \frac{\pi}{L} \frac{1}{2}|\bq|(\bl-\n)\cdot\bsigma}  
          \,\rd \bsigma \,\rd \bq \,\rd \bk\\[6pt]
        =& \sum_{\substack{|\bl|,|\m|, |\n|= -\frac{N}{2},\\ \bl+\n=\bj-\m}}^{\frac{N}{2}-1}  G_3(\bl,\n)\hat{f}_{\bl} \hat{f}_{\m} \hat{f}_{\n}  \,,
    \end{split}
\end{equation}
where
\begin{equation*}
    G_3(\bl,\n):= \int_{\B_R} \int_{\mathbb{S}^{d-1}} |\bq|^{d-2}
        \e^{-\im \frac{\pi}{L} \frac{1}{2} (\bl+\n) \cdot \bq} 
        \e^{\im \frac{\pi}{L} \frac{1}{2}|\bq| (\bl-\n)\cdot\bsigma} \,\rd \bsigma \,\rd \bq.
\end{equation*}

For $\hat{\mathcal{K}}_{4,\bj}^R$, substituting \eqref{fj} into $\hat{\mathcal{K}}_{4,\bj}$ gives:
\begin{equation}\label{KK4}
    \begin{split}
        \hat{\mathcal{K}}_{4,\bj}^R 
        =& \frac{1}{(2L)^{d}} \int_{\D_L} \int_{\B_R} \int_{\mathbb{S}^{d-1}} |\bq|^{d-2}
        \left( \sum_{|\bl| = -\frac{N}{2}}^{\frac{N}{2}-1} \hat{f}_{\bl} \e^{\im \frac{\pi}{L} \bl \cdot \bk'} \right)
        \left(\sum_{|\m| = -\frac{N}{2}}^{\frac{N}{2}-1} \hat{f}_{\m} \e^{\im \frac{\pi}{L} \m \cdot (\bk-\bq)}\right)\\
        & \times 
        \left(\sum_{|\n| = -\frac{N}{2}}^{\frac{N}{2}-1} \hat{f}_{\n} \e^{\im \frac{\pi}{L} \n \cdot \bk}\right)  \e^{-\im \frac{\pi}{L} \bj \cdot \bk} \,\rd \bsigma \,\rd \bq  \,\rd \bk\\[6pt]
        =& \frac{1}{(2L)^{d}} \sum_{|\bl| = -\frac{N}{2}}^{\frac{N}{2}-1} \sum_{|\m| = -\frac{N}{2}}^{\frac{N}{2}-1} \sum_{|\n| = -\frac{N}{2}}^{\frac{N}{2}-1} \hat{f}_{\bl} \hat{f}_{\m} \hat{f}_{\n} \int_{\D_L} \int_{\B_R} \int_{\mathbb{S}^{d-1}} |\bq|^{d-2}\\
        &\times   \e^{\im \frac{\pi}{L} (\bl+\m+\n-\bj) \cdot \bk}
          \e^{-\im \frac{\pi}{L} \frac{1}{2} (2\m+\bl) \cdot \bq} 
          \e^{\im \frac{\pi}{L} \frac{1}{2}|\bq|\bl\cdot\bsigma}  
          \,\rd \bsigma \,\rd \bq \,\rd \bk\\[6pt]
        =& \sum_{\substack{|\bl|,|\m|, |\n|= -\frac{N}{2},\\ \bl+\n=\bj-\m}}^{\frac{N}{2}-1}  G_4(\bl,\m)\hat{f}_{\bl} \hat{f}_{\m} \hat{f}_{\n} \,,
    \end{split}
\end{equation}
where 
\begin{equation*}
    G_4(\bl,\m):= \int_{\B_R} \int_{\mathbb{S}^{d-1}} |\bq|^{d-2}
        \e^{-\im \frac{\pi}{L} \frac{1}{2} (2\m+\bl) \cdot \bq} 
        \e^{\im \frac{\pi}{L} \frac{1}{2}|\bq| \bl\cdot\bsigma} \,\rd \bsigma \,\rd \bq.
\end{equation*}

\subsection{Fast algorithm}
\label{subsec:fast}
Looking back, the most costly component in evaluating $\hat {\mathcal K}_{1,\bj}^R$ in \eqref{KK1} has a complexity of $\mathcal O(N^{3d})$, due to the triple summation that appears in the expression. To reduce this cost, we aim to extract a convolutional structure. 
The central idea of the present method is to approximate the weight functions that includes oscillatory collision kernel by a low-rank separated representation. This is achieved by introducing numerical quadratures in the radial and angular variables appearing in the Boltzmann-type representation of the collision operator.
Specifically, this requires a suitable decomposition for the weight function $G_1(\bl,\m,\n)$, and we seek
\begin{equation} \label{315}
    G_1(\bl,\m,\n) \approx \sum_{p=1}^M \alpha_p(\m) \beta_p(\bl) \gamma_p(\n) \,.
\end{equation}
If this decomposition is successful, \eqref{KK1} can then be rewritten as a summation with a double-convolutional structure:
\begin{equation*}
    \hat{\mathcal{K}}_{1,\bj}^R = \sum_{p=1}^M \sum_{|\m|= -\frac{N}{2}}^{\frac{N}{2}-1}  \left( \alpha_p(\m) \hat{f}_{\m} \right) \underbrace{\sum_{\substack{|\bl|, |\n|= -\frac{N}{2},\\ \bl+\n=\bj-\m}}^{\frac{N}{2}-1} \left( \beta_p(\bl)\hat{f}_{\bl} \right) \left( \gamma_p(\n)\hat{f}_{\n}\right)}_{:=F_p(\bj-\m)}\,.
\end{equation*}
As written, $F_p(\bj-\m)$ is a convolution of  $\beta_p(\bl)\hat{f}_{\bl}$ and $\gamma_p(\n)\hat{f}_{\n}$, and once this is obtained, along with $\alpha_p(\m) \hat{f}_{\m}$,  the summation over $\m$ becomes a convolution as well. As a result, the total cost of evaluating $\hat{\mathcal{K}}_{1,\bj}^{R}$ (for all $\bj$) is reduced from $O(N^{3d})$ to $O(M N^d \log N), M \ll N^{2d-1}$ with the help of FFTs.

To achieve \eqref{315}, we decompose $\bq \in \mathcal{B}_R$ in $|\bq| \in [0,R]$ and $\hat{\bq} \in \mathbb{S}^{d-1}$ with a radial transformation, and approximate the integrals using the quadrature rule:
Let $N_r$ denote the number of quadrature points in the radial direction, $N_s$ the number of quadrature points for the spherical direction, and $N_{sig}$ the number of quadrature points for the collision angle.
\begin{equation}\label{G1}
    \begin{split}
        &G_1(\bl,\m,\n) \\
        =&\int_{\mathcal{B}_R} |\bq|^{d-2} \e^{-\im \frac{\pi}{L} \frac{1}{2} (2\m+\bl+\n) \cdot \bq}
        \left[\int_{\mathbb{S}^{d-1}}  
        \e^{\im \frac{\pi}{L} \frac{1}{2}|\bq|(\bl-\n)\cdot\bsigma}
         \,\rd \bsigma \right] 
         \,\rd \bq \\[6pt]
         =& \int_0^R  |\bq|^{d-2}  \left[ \int_{\mathbb{S}^{d-1}}
         \underbrace{\e^{-\im \frac{\pi}{L} |\bq| \m \cdot \hat{\bq}}}_{:=\alpha(\m,|\bq|,\hat{\bq})} 
        \underbrace{\e^{-\im \frac{\pi}{L} \frac{1}{2} |\bq| \bl \cdot \hat{\bq}}}_{:=\beta(\bl,|\bq|,\hat{\bq})} 
        \underbrace{\e^{-\im \frac{\pi}{L} \frac{1}{2} |\bq| \n \cdot \hat{\bq}}}_{:=\gamma(\n,|\bq|,\hat{\bq})} \,\rd \hat{\bq} \right]\\
        & \qquad \qquad \qquad \qquad \qquad \qquad \times \left[\int_{\mathbb{S}^{d-1}} 
        \underbrace{\e^{\im \frac{\pi}{L} \frac{1}{2}|\bq| \bl \cdot \bsigma}}_{:=\tilde{\beta}(\bl,|\bq|,\bsigma)} 
        \underbrace{\e^{-\im \frac{\pi}{L} \frac{1}{2}|\bq| \n \cdot \bsigma}}_{:=\tilde{\gamma}(\n,|\bq|,\bsigma)}
         \,\rd \bsigma \right] |\bq|^{d-1} \,\rd |\bq|\\[6pt]
         \approx & \sum_{p_1}^{N_{r}}  w_{r}   |\bq|_{p_1}^{2d-3}\Bigg[\sum_{p_2}^{N_{s}}w_{s} \alpha_{p_1p_2}(\m) \beta_{p_1p_2}(\bl) \gamma_{p_1p_2}(\n)\Bigg] \Bigg[\sum_{p_3}^{N_{sig}}w_{sig} \tilde{\beta}_{p_1p_3}(\bl) \tilde{\gamma}_{p_1p_3}(\n)\Bigg],
    \end{split}
\end{equation}
where $w_r$, $w_s$ and $w_{sig}$ are the corresponding quadrature weights for each integral in the radial part $|\bq|$, spherical part $\hat{\bq}$ and $\bsigma$-part. In practice, we use Gauss-Legendre quadrature with $N_r \leq N$ to discretize $|\bq|$, and utilize the mid-point quadrature (in 2D) and a spherical quadrature called Spherical Design \cite{Womersley} (in 3D) with $N_s, N_{sig} \ll {N^{d-1}}$ to discretize $\hat{\bq}$ and $\bsigma$. 

Putting everything together, $\hat{\mathcal{K}}_{1,\bj}^R$ in \eqref{KK1} can be calculated using the following double-convolutional structure:
\begin{equation}\label{K1_con}
    \begin{split}
        \hat{\mathcal{K}}_{1,\bj}^R =& \sum_{\substack{|\bl|,|\m|,|\n| = -\frac{N}{2},\\ \bl+\n=\bj-\m}}^{\frac{N}{2}-1}  G_1(\bl,\m,\n) \hat{f}_{\bl} \hat{f}_{\m} \hat{f}_{\n} \\[6pt]
        =& \sum_{p_1}^{N_{r}}  w_{r} |\bq|_{p_1}^{2d-3} \sum_{p_2}^{N_{s}}w_{s}  \sum_{\substack{|\m|= -\frac{N}{2}}}^{\frac{N}{2}-1} \left( \alpha_{p_1p_2}(\m)\hat{f}_{\m} \right)  \\
        & \qquad \quad \times 
        \underbrace{\Bigg[ \sum_{p_3}^{N_{sig}}w_{sig}
        \sum_{\substack{|\bl|,|\n|= -\frac{N}{2},\\ \bl+\n=\bj-\m}}^{\frac{N}{2}-1} \left( \beta_{p_1p_2}(\bl) \tilde{\beta}_{p_1p_3}(\bl) \hat{f}_{\bl} \right)
        \left( \gamma_{p_1p_2}(\n) \tilde{\gamma}_{p_1p_3}(\n) \hat{f}_{\n} \right) \Bigg]}_{:=F_1(\bj-\m)}.
    \end{split}
\end{equation}


Similar to $G_1$, we seek a decomposition for $G_2$:
\begin{equation*} 
    \begin{split}
        G_2(\m,\n) 
         =& \int_0^R  |\bq|^{d-2} \left[ \int_{\mathbb{S}^{d-1}}
         \underbrace{\e^{-\im \frac{\pi}{L} |\bq| \m \cdot \hat{\bq}}}_{:=\alpha(\m,|\bq|,\hat{\bq})} 
        \underbrace{\e^{-\im \frac{\pi}{L} \frac{1}{2} |\bq| \n \cdot \hat{\bq}}}_{:=\gamma(\n,|\bq|,\hat{\bq})} \,\rd \hat{\bq} \right]\\
        & \qquad \qquad \qquad \qquad \qquad \times 
        \left[\int_{\mathbb{S}^{d-1}} 
        \underbrace{\e^{-\im \frac{\pi}{L} \frac{1}{2}|\bq| \n \cdot \bsigma}}_{:=\tilde{\gamma}(\n,|\bq|,\bsigma)}
         \,\rd \bsigma \right] |\bq|^{d-1} \,\rd |\bq|\\[6pt]
         \approx & \sum_{p_1}^{N_{r}}  w_{r}   |\bq|_{p_1}^{2d-3}
         \Bigg[\sum_{p_2}^{N_{s}}w_{s} \alpha_{p_1p_2}(\m) \gamma_{p_1p_2}(\n)\Bigg] \Bigg[\sum_{p_3}^{N_{sig}}w_{sig} \tilde{\gamma}_{p_1p_3}(\n)\Bigg]\,.
    \end{split}
\end{equation*}
Then
$\hat{\mathcal{K}}_{2,\bj}^R$ in \eqref{KK2} can be calculated in the following convolutional structure:
\begin{multline*}
    \hat{\mathcal{K}}_{2,\bj}^R
        = \sum_{p_1}^{N_{r}}  w_{r} |\bq|_{p_1}^{2d-3} \sum_{p_2}^{N_{s}}w_{s}  \sum_{\substack{|\m|= -\frac{N}{2}}}^{\frac{N}{2}-1}    \left( \alpha_{p_1p_2}(\m)\hat{f}_{\m} \right) \\
        \times 
        \underbrace{\Bigg[ \sum_{p_3}^{N_{sig}}w_{sig}
        \sum_{\substack{|\bl|,|\n|= -\frac{N}{2},\\ \bl+\n=\bj-\m}}^{\frac{N}{2}-1} \hat{f}_{\bl}
        \left( \gamma_{p_1p_2}(\n) \tilde{\gamma}_{p_1p_3}(\n) \hat{f}_{\n} \right) \Bigg]}_{:=F_2(\bj-\m)}\,.
\end{multline*}

Likewise,  $G_3$ can be decomposed as:
\begin{equation*} 
    \begin{split}
        G_3(\bl,\n) =& \int_0^R  |\bq|^{d-2}  \left[ \int_{\mathbb{S}^{d-1}}
         \underbrace{\e^{-\im \frac{\pi}{L} \frac{1}{2} |\bq| \bl \cdot \hat{\bq}}}_{:=\beta(\bl,|\bq|,\hat{\bq})} 
        \underbrace{\e^{-\im \frac{\pi}{L} \frac{1}{2} |\bq| \n \cdot \hat{\bq}}}_{:=\gamma(\n,|\bq|,\hat{\bq})} \,\rd \hat{\bq} \right]\\
        &\qquad \qquad \qquad \qquad \qquad \quad \times 
        \left[\int_{\mathbb{S}^{d-1}} 
        \underbrace{\e^{\im \frac{\pi}{L} \frac{1}{2}|\bq| \bl \cdot \bsigma}}_{:=\tilde{\beta}(\bl,|\bq|,\bsigma)}
        \underbrace{\e^{-\im \frac{\pi}{L} \frac{1}{2}|\bq| \n \cdot \bsigma}}_{:=\tilde{\gamma}(\n,|\bq|,\bsigma)}
         \,\rd \bsigma \right] |\bq|^{d-1} \,\rd |\bq|\\[6pt]
         \approx & \sum_{p_1}^{N_{r}}  w_{r}  |\bq|_{p_1}^{2d-3}\Bigg[\sum_{p_2}^{N_{s}}w_{s} \beta{p_1p_2}(\bl) \gamma_{p_1p_2}(\n)\Bigg] \Bigg[\sum_{p_3}^{N_{sig}}w_{sig} \tilde{\beta}_{p_1p_3}(\bl)\tilde{\gamma}_{p_1p_3}(\n)\Bigg].
    \end{split}
\end{equation*}
Thus, the computation of 
$\hat{\mathcal{K}}_{3,\bj}^R$ in \eqref{KK3} can be accelerated by solving the following convolutional structure via FFT:
\begin{multline*}
    \hat{\mathcal{K}}_{3,\bj}^R =\sum_{p_1}^{N_{r}}  w_{r} B_{p_1} |\bq|_{p_1}^{d-1} \sum_{p_2}^{N_{s}}w_{s}  \sum_{\substack{|\m|= -\frac{N}{2}}}^{\frac{N}{2}-1}   \hat{f}_{\m} \\
    \times 
        \underbrace{\Bigg[ \sum_{p_3}^{N_{sig}}w_{sig}
        \sum_{\substack{|\bl|,|\n|= -\frac{N}{2},\\ \bl+\n=\bj-\m}}^{\frac{N}{2}-1} \left( \beta{p_1p_2}(\bl) \tilde{\beta}_{p_1p_3}(\bl)\hat{f}_{\bl} \right)
        \left( \gamma_{p_1p_2}(\n) \tilde{\gamma}_{p_1p_3}(\n) \hat{f}_{\n} \right) \Bigg]}_{:=F_3(\bj-\m)}.
\end{multline*}

Similarly for $G_4$, we have the following decomposition: 
\begin{equation*} 
    \begin{split}
        G_4(\bl,\m) =& \int_0^R  |\bq|^{d-2}  \left[ \int_{\mathbb{S}^{d-1}}
         \underbrace{\e^{-\im \frac{\pi}{L} |\bq| \m \cdot \hat{\bq}}}_{:=\alpha(\m,|\bq|,\hat{\bq})} 
        \underbrace{\e^{-\im \frac{\pi}{L} \frac{1}{2} |\bq| \bl \cdot \hat{\bq}}}_{:=\beta(\bl,|\bq|,\hat{\bq})} 
        \,\rd \hat{\bq} \right]\\
        & \qquad \qquad \qquad \qquad\qquad \qquad \qquad  \times
        \left[\int_{\mathbb{S}^{d-1}} 
        \underbrace{\e^{\im \frac{\pi}{L} \frac{1}{2}|\bq| \bl \cdot \bsigma}}_{:=\tilde{\beta}(\bl,|\bq|,\bsigma)} 
         \,\rd \bsigma \right] |\bq|^{d-1} \,\rd |\bq|\\[6pt]
         \approx & \sum_{p_1}^{N_{r}}  w_{r}   |\bq|_{p_1}^{2d-3}\Bigg[\sum_{p_2}^{N_{s}}w_{s} \alpha_{p_1p_2}(\m) \beta_{p_1p_2}(\bl)\Bigg] \Bigg[\sum_{p_3}^{N_{sig}}w_{sig} \tilde{\beta}_{p_1p_3}(\bl)\Bigg],
    \end{split}
\end{equation*}
and $\hat{\mathcal{K}}_{4,\bj}^R$ in \eqref{KK4} has the  convolutional structure, which can be evaluated by FFT: 
\begin{multline*}
    \hat{\mathcal{K}}_{4,\bj}^R 
        = \sum_{p_1}^{N_{r}}  w_{r} B_{p_1} |\bq|_{p_1}^{d-1} \sum_{p_2}^{N_{s}}w_{s}  \sum_{\substack{|\m|= -\frac{N}{2}}}^{\frac{N}{2}-1}    \left( \alpha_{p_1p_2}(\m)\hat{f}_{\m} \right)\\  
        \underbrace{\Bigg[ \sum_{p_3}^{N_{sig}}w_{sig}
        \sum_{\substack{|\bl|,|\n|= -\frac{N}{2},\\ \bl+\n=\bj-\m}}^{\frac{N}{2}-1} \left( \beta_{p_1p_2}(\bl) \tilde{\beta}_{p_1p_3}(\bl) \hat{f}_{\bl} \right)
        \hat{f}_{\n} \Bigg]}_{:=F_4(\bj-\m)}\,.
\end{multline*}


\subsection{Implementation and complexity analysis}
\label{sub:complexity}

To better clarify our fast algorithm, the pseudo code is presented as in Algorithm \ref{alg:fast-wke}, which summarizes the implementation of the fast Fourier spectral evaluation of the collision operator. The algorithm clearly separates the Fourier
discretization parameter $N$ from the kernel-approximation parameters $(N_r,N_s,N_{sig})$, which play distinct roles in accuracy and efficiency.
\begin{algorithm}[htbp]
\caption{Fast Fourier spectral evaluation of the collision operator of 4-WKE}
\label{alg:fast-wke}
\begin{algorithmic}[1]
\REQUIRE
\STATE $f(\bk)$: The distribution function (array of size $N^d$) on a uniform grid in $\mathcal{D}_L$;
\STATE $N_r, N_s, N_{sig}$: Number of quadrature points for radial part $|\bq|$, spherical part $\hat{\bq}$, and collision angle $\bsigma$ integrations;
\STATE $\{|\bq|_{p_1}, w_r\}_{p_1=1}^{N_r}$,
$\{\hat \bq_{p_2}, w_s\}_{p_2=1}^{N_s}$,
$\{\bsigma_{p_3}, w_{sig}\}_{p_3=1}^{N_{sig}}$: Corresponding quadrature nodes and weights;
\STATE $\alpha, \beta, \gamma$: Pre-computed phase factors from kernel decomposition.
\ENSURE Collision operator $\mathcal{K}^R(f,f,f)$ on the $\bk$-grid.

\STATE Compute Fourier modes $\hat f = \hat f_\bj$ by FFT for $- \frac{N}{2} \leq |\bj| \leq \frac{N}{2}-1$.

\STATE Initialize $\mathcal{\hat K}_{1,\bj}^R = \mathcal{\hat K}_{2,\bj}^R = \mathcal{\hat K}_{3,\bj}^R = \mathcal{\hat K}_{4,\bj}^R = 0$.

\FOR{$p_1=1,\dots,N_r$}
    \STATE $|\bq| \leftarrow |\bq|_{p_1}$
    \FOR{$p_2=1,\dots,N_s$}
        \STATE $\hat \bq \leftarrow \hat \bq_{p_2}$
        \STATE Evaluate phase factors in the Fourier space:
        \[
        \alpha(\m,|\bq|,\hat{\bq}) = \e^{-\im \frac{\pi}{L} |\bq| \m \cdot \hat{\bq}} \leftarrow \alpha_{p_1p_2}(\m)
        \]
        \[
        \beta(\bl,|\bq|,\hat{\bq}) = \e^{-\im \frac{\pi}{L} \frac{1}{2} |\bq| \bl \cdot \hat{\bq}}  \leftarrow \beta_{p_1p_2}(\bl)
        \]
        \[
        \gamma(\n,|\bq|,\hat{\bq}) = \e^{-\im \frac{\pi}{L} \frac{1}{2} |\bq| \n \cdot \hat{\bq}} \leftarrow \gamma_{p_1p_2}(\n)
        \]
        \STATE Form weighted spectra:
        \[
        \hat F_{\m} = \alpha_{p_1p_2}(\m) \hat f,\quad
        \hat F_{\bl} = \beta_{p_1p_2}(\bl) \hat f,\quad
        \hat F_{\n} = \gamma_{p_1p_2}(\n) \hat f
        \]
        \FOR{$p_3=1,\dots,N_{sig}$}
            \STATE $\bsigma \leftarrow \bsigma_{p_3}$
            \STATE Evaluate additional phase factors in Fourier space:
            \[
            \tilde{\beta}(\bl,|\bq|,\bsigma) = \e^{\im \frac{\pi}{L} \frac{1}{2}|\bq| \bl \cdot \bsigma} \leftarrow \tilde{\beta}_{p_1p_3}(\bl)
            \]
            \[
            \tilde{\gamma}(\n,|\bq|,\hat{\bq}) = \e^{-\im \frac{\pi}{L} \frac{1}{2} |\bq| \n \cdot \hat{\bq}} \leftarrow \tilde{\gamma}_{p_1p_3}(\n)
            \]
            \STATE Form weighted spectra:
            \[
            \hat G_{\bl} = [\beta_{p_1p_2}(\bl) \tilde{\beta}_{p_1p_3}(\bl)] \hat f, \quad
            \hat G_{\n} = [\gamma_{p_1p_2}(\n)  \tilde{\gamma}_{p_1p_3}(\n)] \hat f
            \]
            \STATE Compute the inverse FFTs:
            \[
            g_{\bl} = \text{iFFT}(\hat G_l),\quad g_{\n} = \text{iFFT}(\hat G_n)
            \]
            \STATE First inner convolution in $(\bl,\n)$:
            \[
            \widehat{(g_{\bl} g_{\n})}=\text{FFT}(g_{\bl} \cdot g_{\n})
            \]
            \STATE Compute the inverse FFT of $\hat F_{\m}$:
            \[
            f_{\m} = \text{iFFT}(\hat F_{\m})
            \]
            \STATE Second outer convolution:
            \[
            \widehat{f_{\m}(g_{\bl} g_{\n})} = \text{FFT}\!\left(f_{\m}\cdot\text{iFFT}\big(\widehat{(g_{\bl} g_{\n})}\big)\right)
            \]
            \STATE Accumulate $\mathcal{\hat K}_{1,\bj}^R$ with quadrature weights according to \eqref{K1_con}, and similarly apply to $\mathcal{\hat K}_{2,\bj}^R, \mathcal{\hat K}_{3,\bj}^R, \mathcal{\hat K}_{4,\bj}^R$.
        \ENDFOR
    \ENDFOR
\ENDFOR

\STATE Combine terms:
\[
\mathcal{\hat K}_\bj^R = \frac{1}{2^{d-1}}(\mathcal{\hat K}_{1,\bj}^R - \mathcal{\hat K}_{2,\bj}^R + \mathcal{\hat K}_{3,\bj}^R - \mathcal{\hat K}_{4,\bj}^R).
\]

\STATE Compute $\mathcal{K}^R(f,f,f)=\text{iFFT}(\mathcal{\hat K}_\bj^R)$.

\RETURN $\mathcal{K}^R(f,f,f)$.
\end{algorithmic}
\end{algorithm}

We conclude this section with an explicit evaluation of the computational cost of the proposed fast algorithm. Let $N$ denote the number of Fourier modes per dimension, so that the total number of degrees of freedom is $N^d$. The low-rank approximation in \eqref{G1} introduces three additional parameters: the number of radial quadrature points $N_r$, spherical quadrature points $N_s$, and collision-angle quadrature points $N_{sig}$.

For each quadrature configuration, the collision operator evaluation consists of a fixed
number of FFT-based convolutions, each with cost $O(N^d\log N)$. Consequently, the overall
computational complexity per time step scales as
\[
O\!\left(N_r N_s N_{sig}\, N^d \log N\right).
\]
In sum, taking the most complicated $\mathcal{K}^R_1$ as example,
\begin{itemize}
    \item No pre-computation and storage requirement for weight functions;
    \item Compute the $\hat{f}_{\bj}$ by using the FFT -- cost $O(N^d \log N)$;
    \item Compute the $\hat{\mathcal{K}}_{1,\bj}^R$ in the double-convolution form \eqref{K1_con} using the FFT -- cost $O(N_r N_s N_{sig} N^d \log N)$ with $N_r \leq N$ and $N_s,N_{sig} \ll N^{d-1}$;
    \item Compute $\mathcal{K}_1^R$ by using the inverse FFT to $\hat{\mathcal{K}}_{1,\bj}^R$ -- cost $O(N^d \log N)$.
\end{itemize}
Note that, in practical implementation, the quadrature parameters $(N_r,N_s,N_{sig})$ are chosen much smaller than $N$ (i.e., $N_r \leq N$ and $N_s, N_{sig} \ll N^{d-1}$) and can often be kept fixed while refining the Fourier resolution as verified in the numerical experiments in Section~\ref{sec:numerical}. This leads to a dramatic reduction in computational cost compared to the $O(N^{3d})$ complexity of a direct spectral evaluation, and makes the proposed method feasible for high-resolution simulations of the 4-WKE.

\begin{remark}
    Note that the approximation \eqref{315} in the fast algorithm introduces a truncation error governed by the quadrature accuracy in the radial, spherical and angular variables. This error can be systematically reduced by increasing the quadrature resolution and follows standard error estimates from quadrature theory \cite{womersley2018efficient}. A rigorous analysis of the resulting fully discrete solution, in particular spectral convergence and long-time accuracy, is however highly nontrivial due to the lack of positivity preservation and entropy structure in Fourier spectral discretizations.  The same issue also exits in  the Fourier spectral method for the Boltzmann equation \cite{FM11, HQY20}, and is left for future work.
\end{remark}

\section{Numerical Tests}
\label{sec:numerical}

In this section, we provide extensive numerical results showcasing the application of the fast Fourier spectral method to solve the WKE \eqref{WKE} in both equilibrium and time-evolving states. The numerical experiments are performed in two dimensions (2D) and three dimensions (3D), respectively. 
\subsection{2D case}
\label{subsec:2D}

\subsubsection{Stationary state test in 2D}
\label{subsubsec:stationary-2D}

We first verify the accuracy of the method by substituting the RJ stationary solution \eqref{equi} and checking the vanishing of the integral.

\begin{figure}[h!]
        \begin{subfigure}[t]{0.495\textwidth}
            \includegraphics[width=\textwidth]{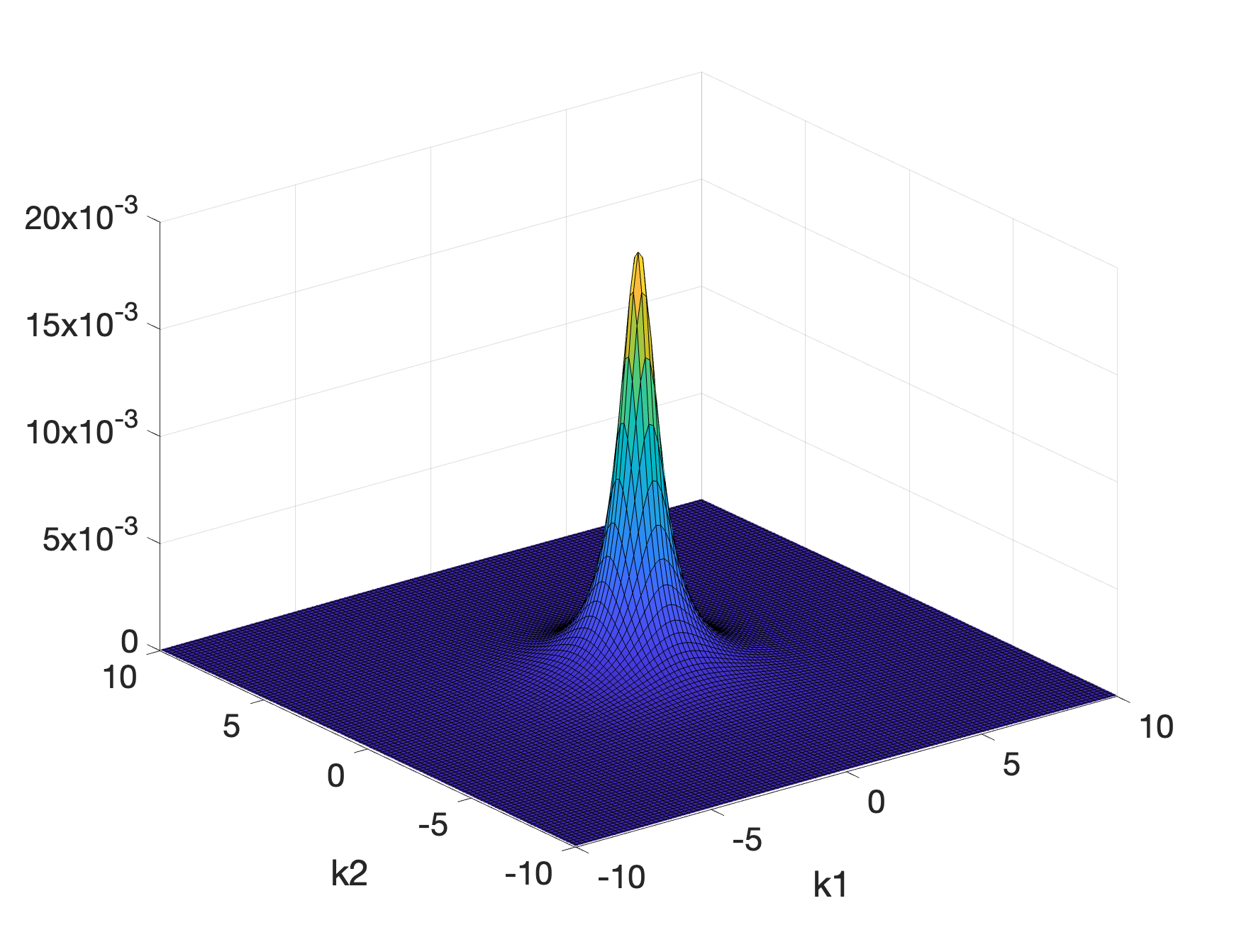}
            \caption{Profile of stationary solution $f_{\text{eq}}$.}
        \end{subfigure}
        \begin{subfigure}[t]{0.495\textwidth}
            \includegraphics[width=\textwidth]{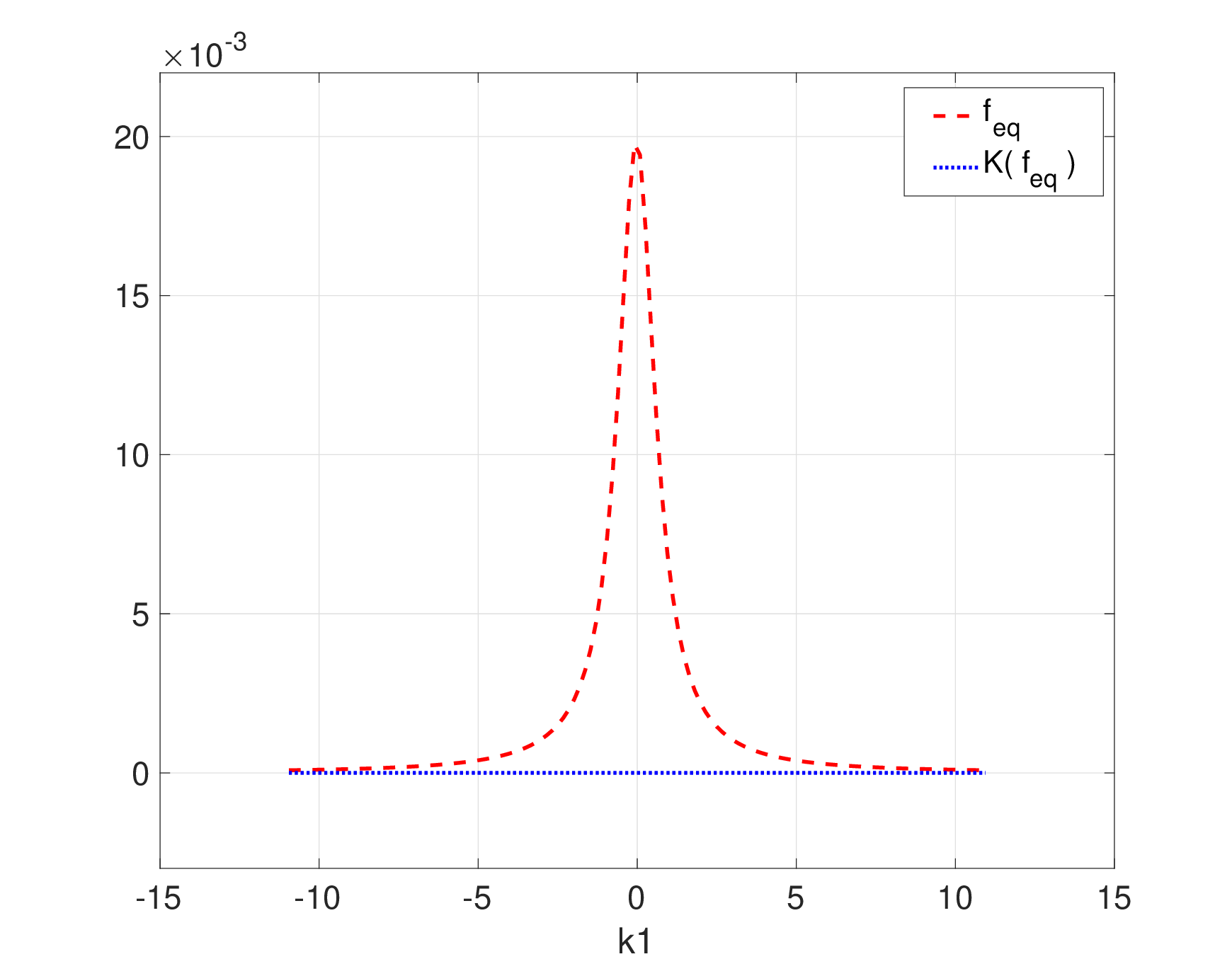}
            \caption{$f_{\text{eq}}$ and $\mathcal{K}(f_{\text{eq}})$.}
        \end{subfigure}
        \caption{Stationary solution in 2D with $N_r = N = 128, N_s = N_{sig} = 12, S = 5$.}
        \label{fig-feq}
\end{figure}

\textit{Example 1} (Stationary solution). 
For the equilibrium solution $f_{\text{eq}}$ in \eqref{equi} with $\mu = 50$, $\bm{\nu} = (5,5)$ and $\xi=100$, we compute $\mathcal{K}(f_{\text{eq}})$ and present the errors in $L^{\infty}$ and $L^{2}$ norms in Table.~\ref{table-2D-Keq}, along with the corresponding computational time in Table.~\ref{table-2D-Keq-2}.
The Fig.~\ref{fig-feq} describes profile of the $f_{\text{eq}}$ and $\mathcal{K}(f_{\text{eq}})$.

\renewcommand{\arraystretch}{1.2}
\begin{table}[h!]
\begin{tabular}{|c|| c | c |}
 \hline
  $N$ & $\| \mathcal{K}(f_{\text{eq}})\|_{L^{\infty}}$ & $\| \mathcal{K}(f_{\text{eq}})\|_{L^{2}}$  \\ [0.5ex]
 \hline
   16 & $1.2126 \times 10^{-7}$ & $9.0791 \times 10^{-8}$ \\  
 \hline
    32 & $9.0791\times 10^{-8}$ & $3.2953\times 10^{-14}$ \\ 
    \hline
    64 & $4.5939\times 10^{-9}$ & $2.8443\times 10^{-17}$ \\ 
    \hline
    128 & $4.4954\times 10^{-10}$ & $9.0718\times 10^{-18}$ \\ 
 \hline
\end{tabular}
\caption{$L^{\infty}$ and $L^2$ error for stationary solution for $N_r = N$ with $ N_s = N_{sig} = 12$ and $S = 5$.}
\label{table-2D-Keq}
\end{table}

\begin{table}[h!]
\begin{tabular}{|c || c | c | c | c || c |} 
 \hline
    $N$ & $\mathcal{K}_1$ & $\mathcal{K}_2$ & $\mathcal{K}_3$ & $\mathcal{K}_4$ &$\mathcal{K}$\\ [0.5ex]
 \hline
 16 & $0.15$s & $0.10$s & $0.12$s & $0.10$s & $0.47$s\\ 
 \hline
 32 & $0.74$s & $0.49$s & $0.58$s & $0.46$s & $2.27$s\\ 
 \hline
  64 & $3.68$s & $2.40$s & $2.86$s & $2.43$s & $11.37$s\\ 
 \hline
  128 & $20.23$s & $13.95$s & $15.92$s & $13.66$s & $63.76$s\\
 \hline
\end{tabular}
\caption{Computational time in 2D for $N_r = N, N_s = N_{sig} = 12 $ and $ S = 5$.}
\label{table-2D-Keq-2}
\end{table}

To illustrate the effect of the parameters $N_s$ and $N_{sig}$, we compute $\mathcal{K}(f_{\text{eq}})$ and show the errors in $L^{\infty}$ and $L^{2}$ norms in Table.~\ref{table-2D-Keq-Nsig} for different $N_s$ and $N_{sig}$, which numerically demonstrate the validity of the approximation via quadrature rules in the fast algorithm introduced in Sec.~\ref{subsec:fast}.

\renewcommand{\arraystretch}{1.2}
\begin{table}[h!]
\begin{tabular}{|c|| c | c |}
 \hline
  $N_s = N_{sig}$ & $\| \mathcal{K}(f_{\text{eq}})\|_{L^{\infty}}$ & $\| \mathcal{K}(f_{\text{eq}})\|_{L^{2}}$  \\ [0.5ex]
 \hline
   3 & $9.3118 \times 10^{-8}$ & $3.1187 \times 10^{-14}$ \\  
 \hline
    6 & $2.6994\times 10^{-8}$ & $5.2885\times 10^{-16}$ \\ 
    \hline
    12 & $4.5939\times 10^{-9}$ & $2.8443\times 10^{-17}$ \\ 
 \hline
\end{tabular}
\caption{$L^{\infty}$ and $L^2$ error for stationary solution in 2D for $N_s = N_{sig}$ with $N = N_{r} = 64$ and $S = 5$.}
\label{table-2D-Keq-Nsig}
\end{table}

\begin{remark}
    Note that rigorously characterizing the exact equilibrium of the semi-discrete scheme is challenging due to the frequency domain truncation and the lack of inherent positivity preservation in Fourier spectral methods. However, the stability of the equilibrium under truncation is supported by recent theoretical results: specifically, in \cite{Menegaki2024}, the solutions initially close to the Rayleigh-Jeans equilibrium $f_{\text{eq}}$ are proved to remain in a close neighborhood in the $L^2$ sense. Furthermore, an alternative route to characterizing the long-time equilibrium is suggested by the analysis of the homogeneous Boltzmann equation in \cite{FM11}, where mass conservation and entropy dissipation lead to convergence toward a constant equilibrium determined by the initial mass. A similar framework could potentially be adapted to the truncated WKE once the relevant conservation and dissipation structures are established, and we leave this to our future work.
\end{remark}

\subsubsection{Time evolution test in 2D}
\label{subsubsec:time-2D}

Now, we perform numerical tests for the 2D equation \eqref{WKE}. The time discretization is performed by the fourth-order Runge-Kutta (RK4) methods.

\begin{figure}[htp]
\begin{adjustwidth}{-1cm}{-1cm}
\centering
        \begin{subfigure}[t]{0.49\linewidth}
            \includegraphics[width=\textwidth]{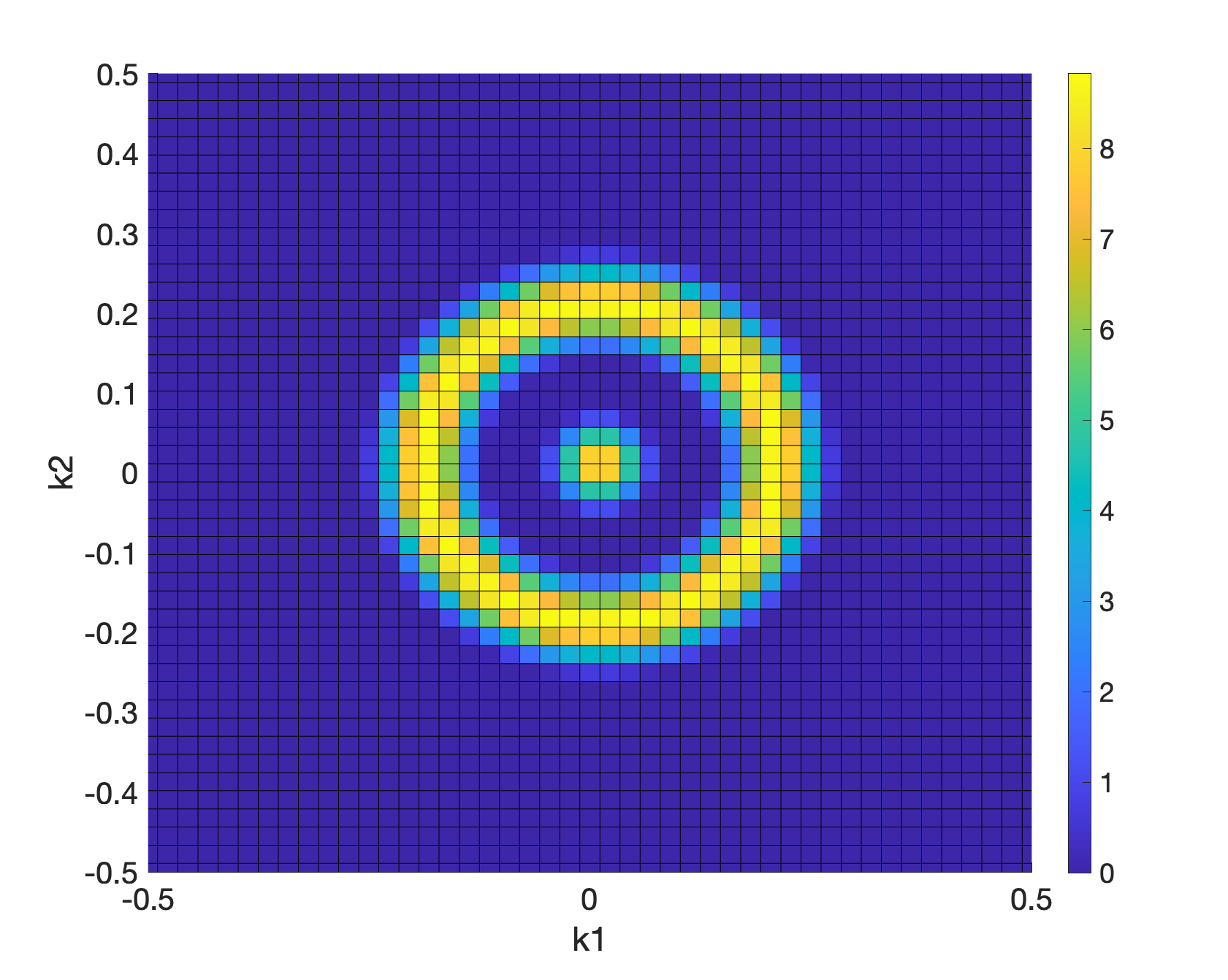}
            \caption{t=0.0}
        \end{subfigure}
        \begin{subfigure}[t]{0.49\linewidth}
            \includegraphics[width=\textwidth]{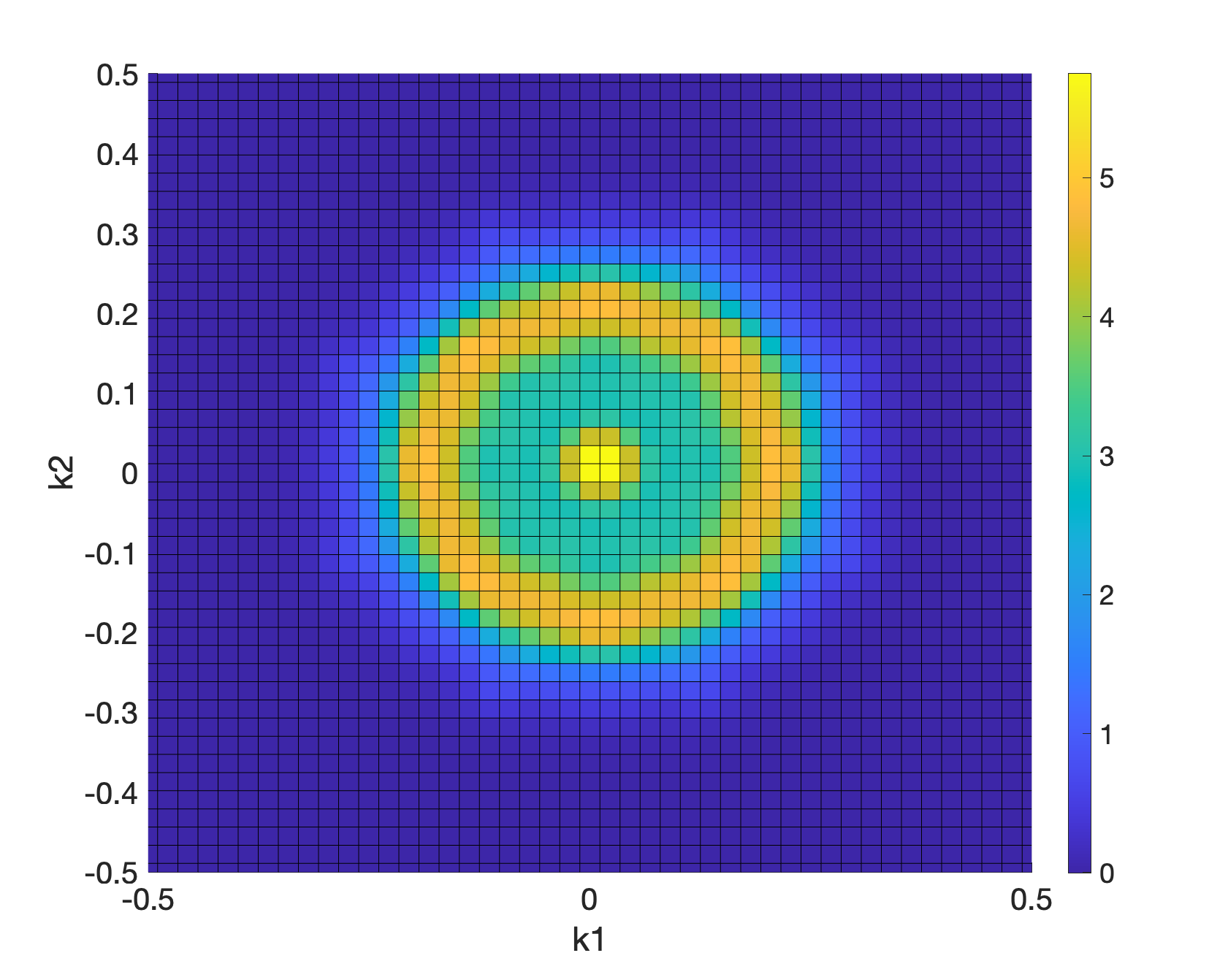}
            \caption{t=0.2}
        \end{subfigure}
        \begin{subfigure}[t]{0.49\linewidth}
            \includegraphics[width=\textwidth]{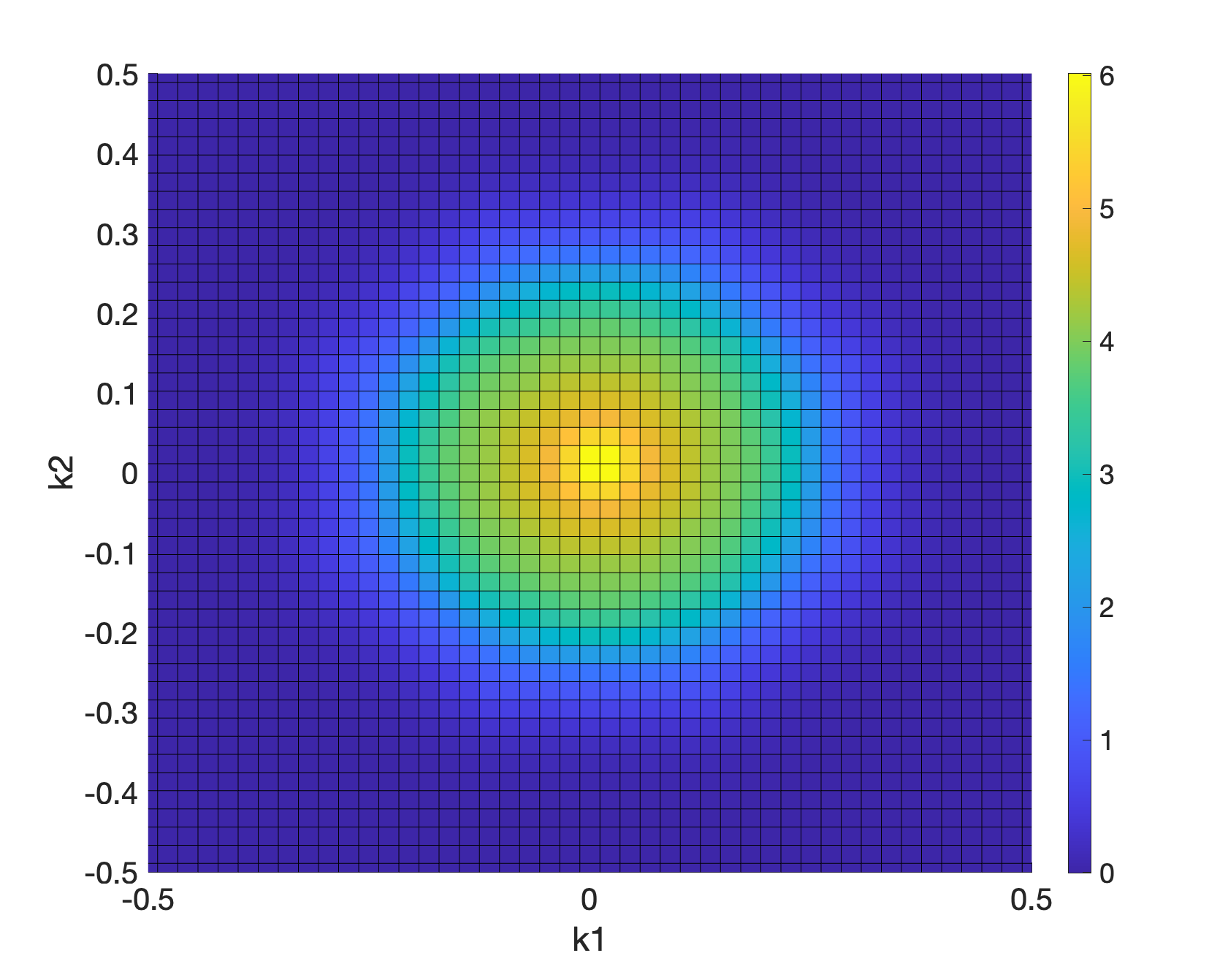}
            \caption{t=0.4}
        \end{subfigure}
        \begin{subfigure}[t]{0.49\linewidth}
            \includegraphics[width=\textwidth]{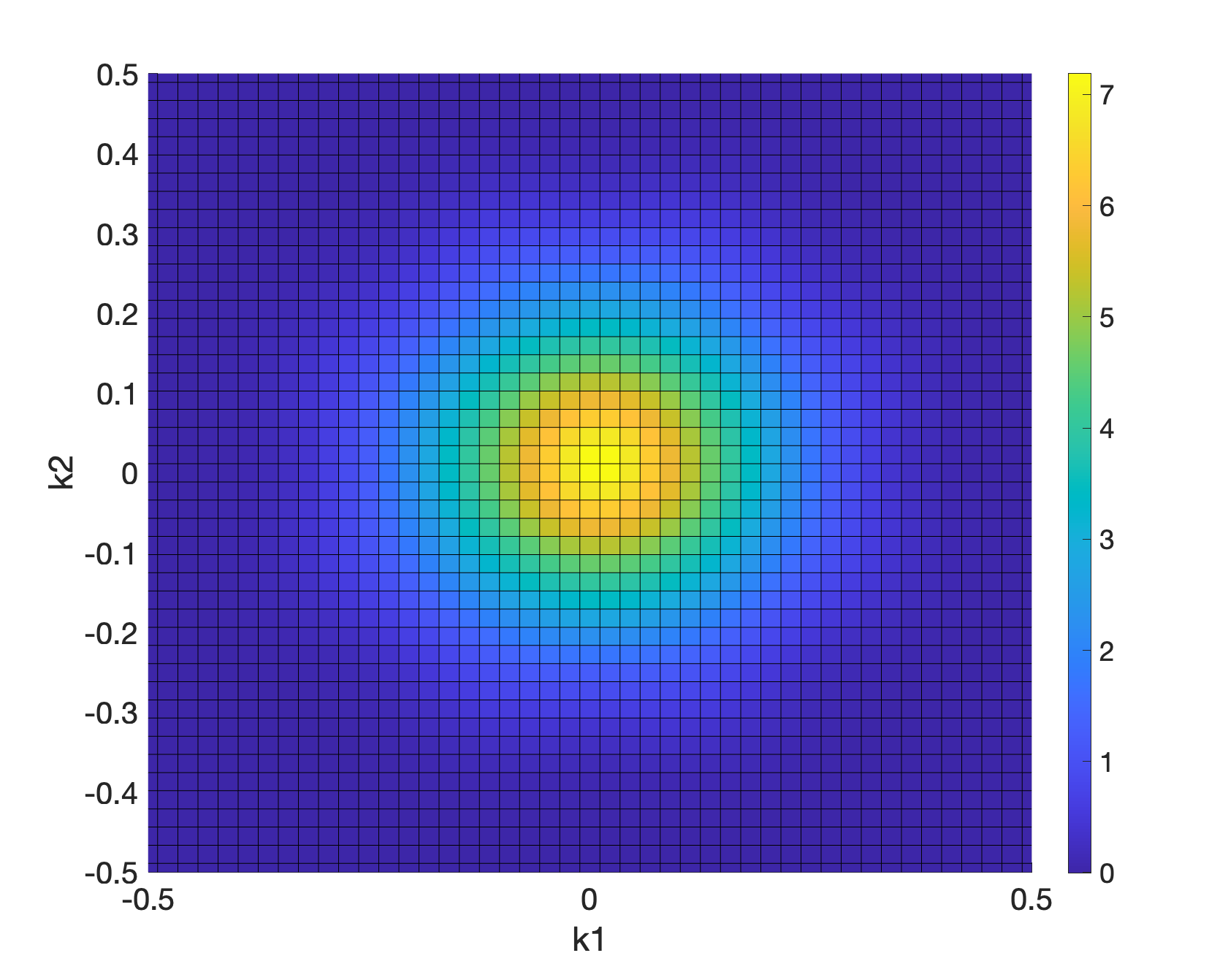}
            \caption{t=0.6}
        \end{subfigure}
        \begin{subfigure}[t]{0.49\linewidth}
            \includegraphics[width=\textwidth]{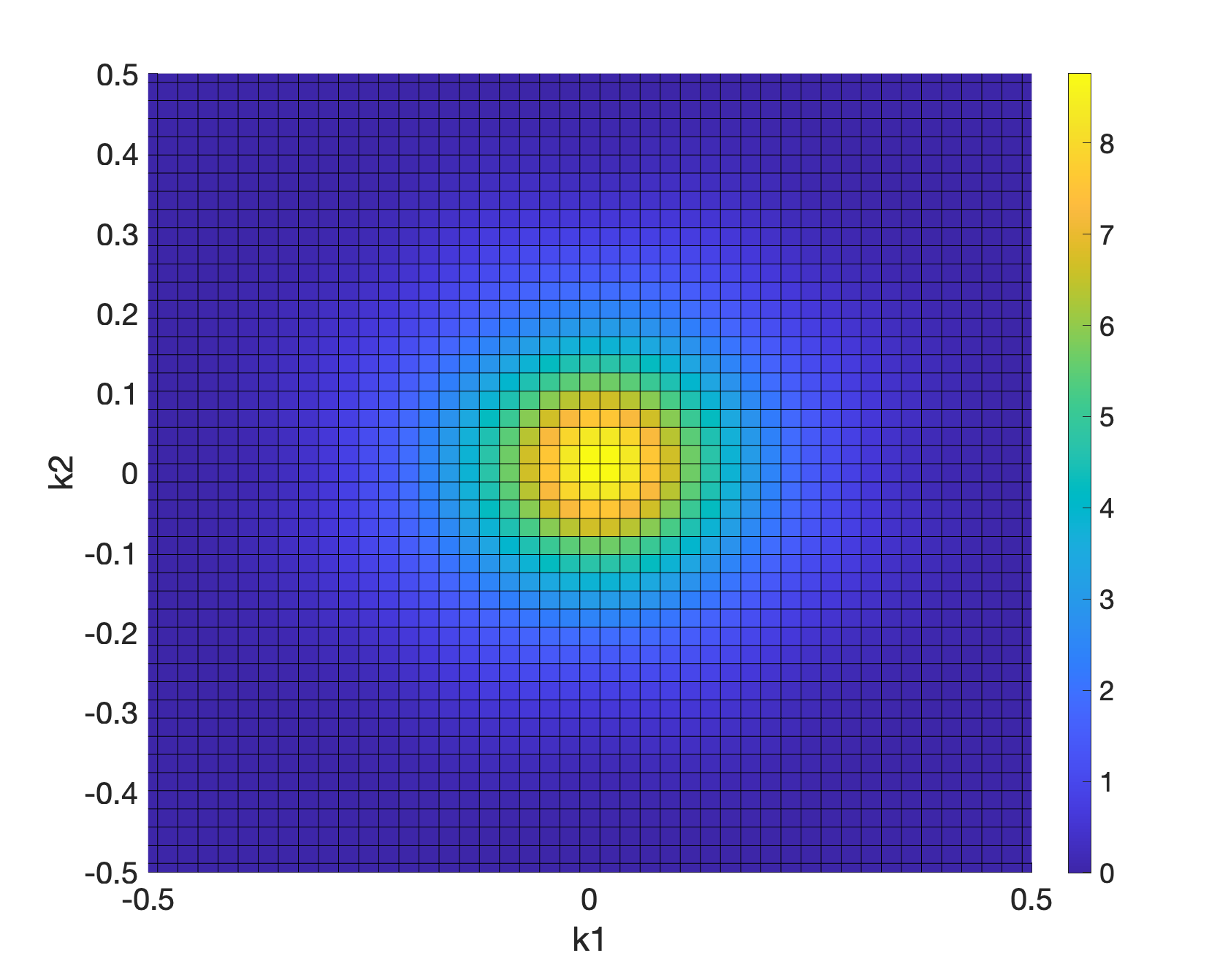}
            \caption{t=0.8}
        \end{subfigure}
        \begin{subfigure}[t]{0.49\linewidth}
            \includegraphics[width=\textwidth]{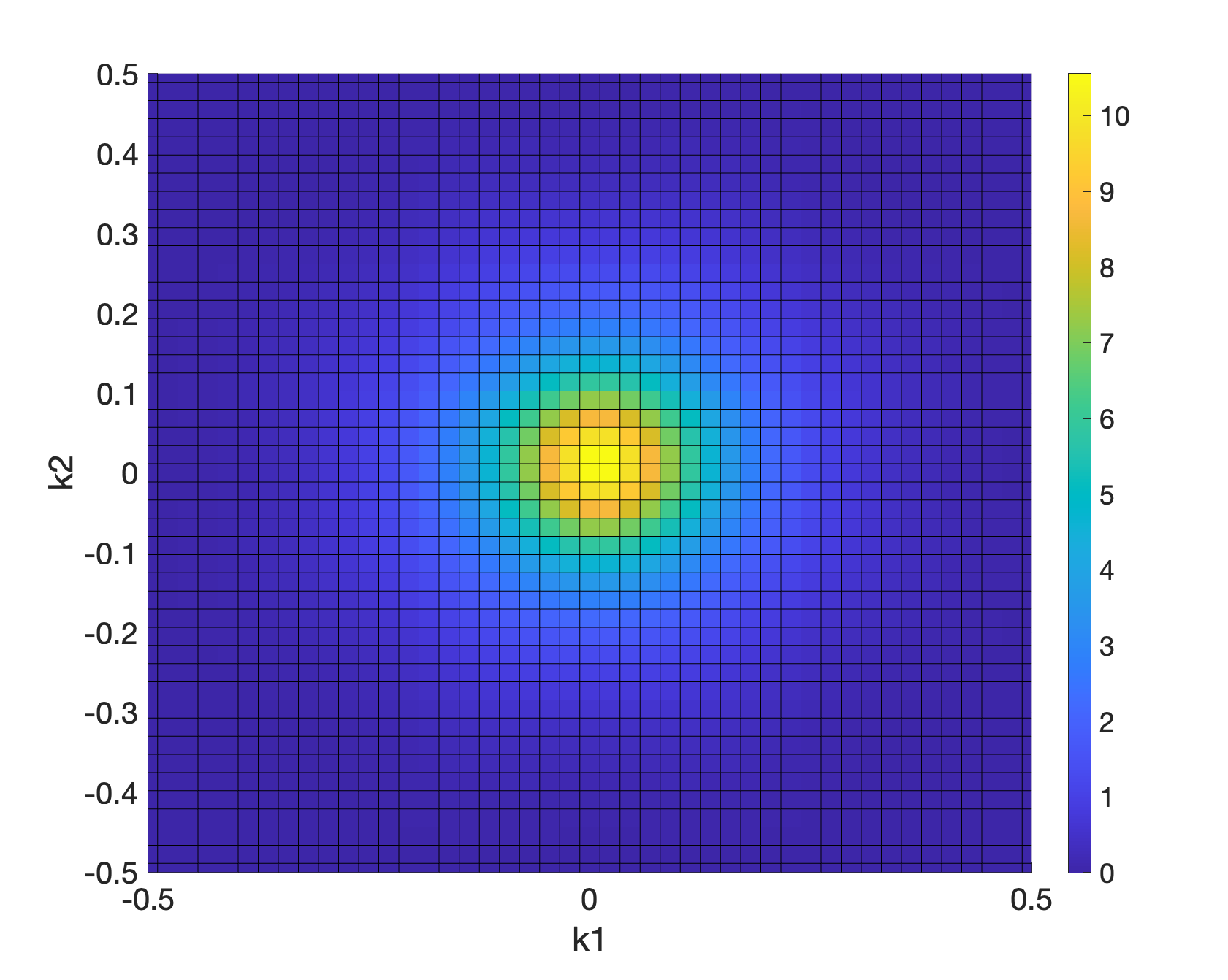}
            \caption{t=1.0}
        \end{subfigure}
\end{adjustwidth}
        \caption{Top-to-bottom view of the time evolution of the solution profile for the isotropic initial condition \eqref{initialdelta} with $N_r = N = 64, N_s = N_{sig} = 12$, $S=0.33$ and $\Delta t = 0.1$.}
        \label{fig-delta-top}
    \end{figure}

\textit{Example 2} (Isotropic case). 
We consider the following  isotropic initial condition:
\begin{equation}\label{initialdelta}
f^0(\bk) = \frac{1}{3} \left( \delta_{u}(|\bk|) + \delta_{u}(|\bk|-0.2)\right),
\end{equation}
where $\delta_{u}(\bk)$ is an approximated delta function given by:
\begin{equation*}
\delta_{u}(\bk)=
\begin{cases}
\frac{1}{2u}\left(1+\cos|\frac{\pi \bk}{u}|\right), \quad & |\bk|\leq u,\\[4pt]
0, & \text{otherwise}\,.
\end{cases}
\end{equation*} 
Here $u$ is taken to be $0.5\sqrt{\Delta \bk}$ with $\Delta \bk$ being the mesh size in $\bk$. 
The time evolution of the solution in the top-to-bottom view is presented in Fig.~\ref{fig-delta-top}. This numerical experiment is consistent with the theoretical result by Escobedo and Vel\'azquez in \cite[Section 1.1.1-1.1.2]{EV2015-1} that the WKE ``generically" blows up, and in the form of a Dirac singularity.

\begin{figure}[h!]
        \begin{subfigure}[t]{0.495\textwidth}
            \includegraphics[width=\textwidth]{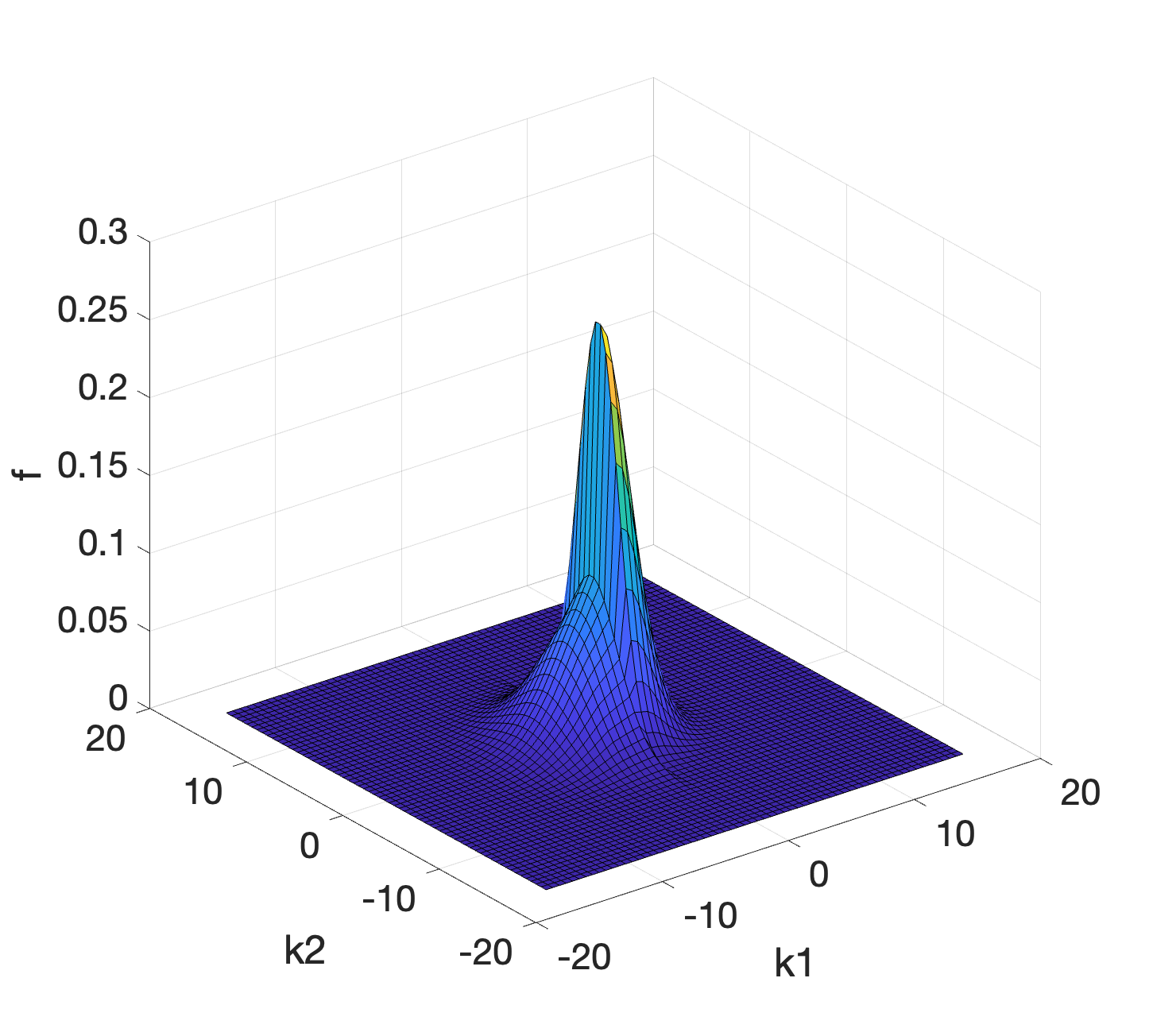}
            \caption{Profile of $f^0_{\text{dis}}$}
        \end{subfigure}
        \begin{subfigure}[t]{0.495\textwidth}
            \includegraphics[width=\textwidth]{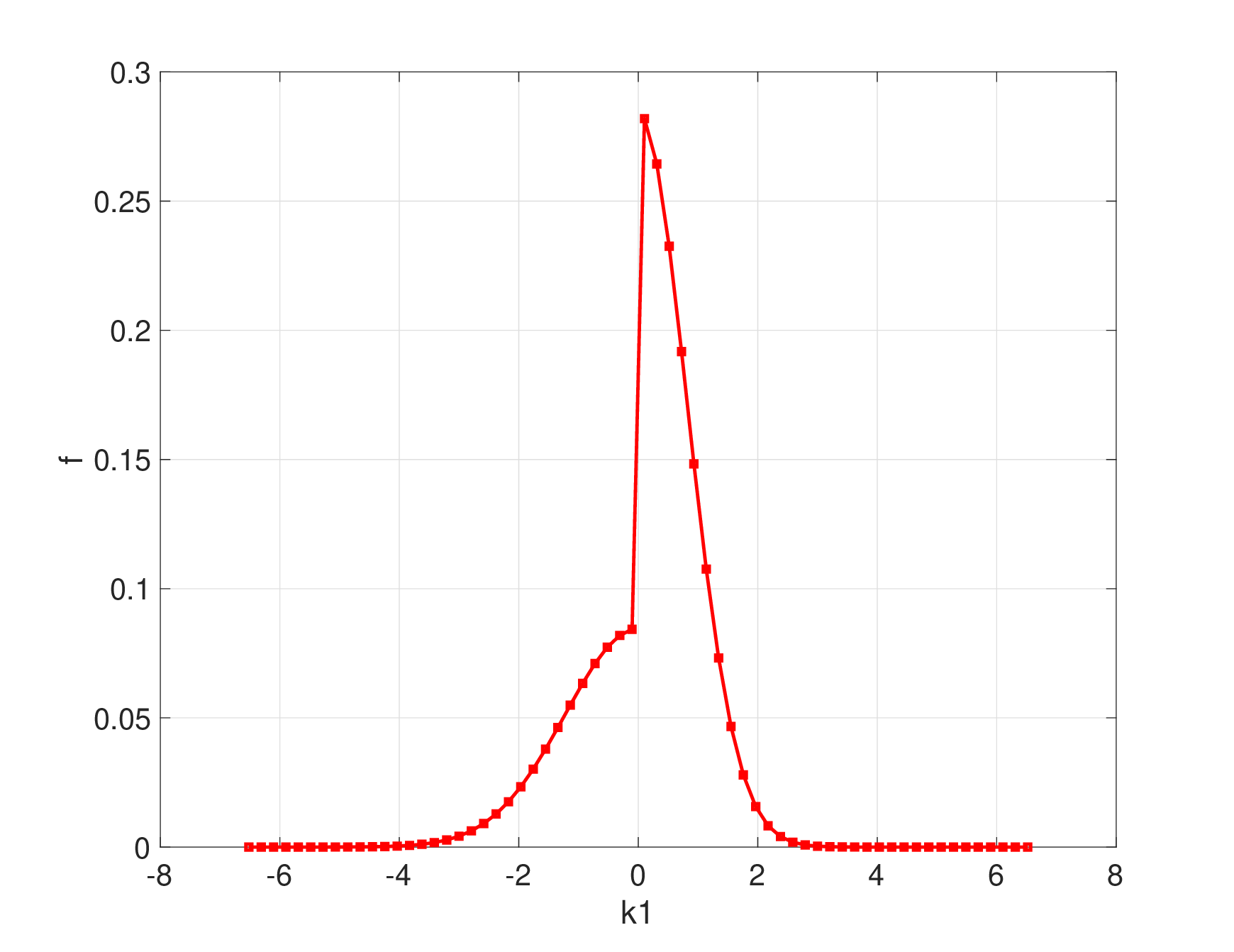}
            \caption{Intersected profile of $f^0_{\text{dis}}(\bk_1, \bk_2 = 0)$}
        \end{subfigure}
        \caption{Discontinous initial condition $f^0_{\text{dis}}$ in 2D for $N = N_r=64, N_s = N_{sig} = 12$ and $S = 3$.}
        \label{fig-disinitial}
    \end{figure}

\begin{figure}[h!]
        \centering
        \begin{subfigure}[t]{0.85\textwidth}
            \includegraphics[width=\textwidth]{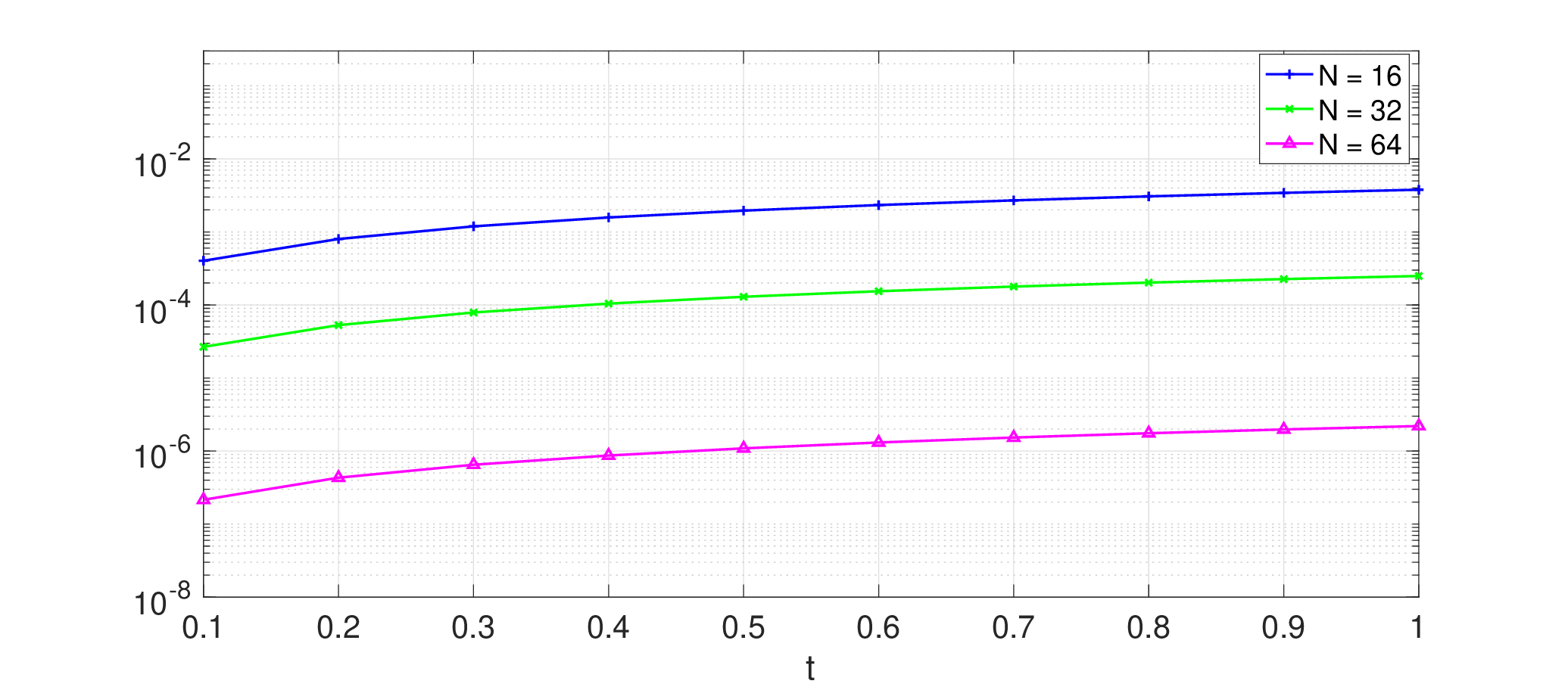}
            \caption{Mass conservation $|\rho(t) - \rho(0)|$ in 2D.}
            \label{fig-mass-2D}
        \end{subfigure}
        \begin{subfigure}[t]{0.85\textwidth}
            \includegraphics[width=\textwidth]{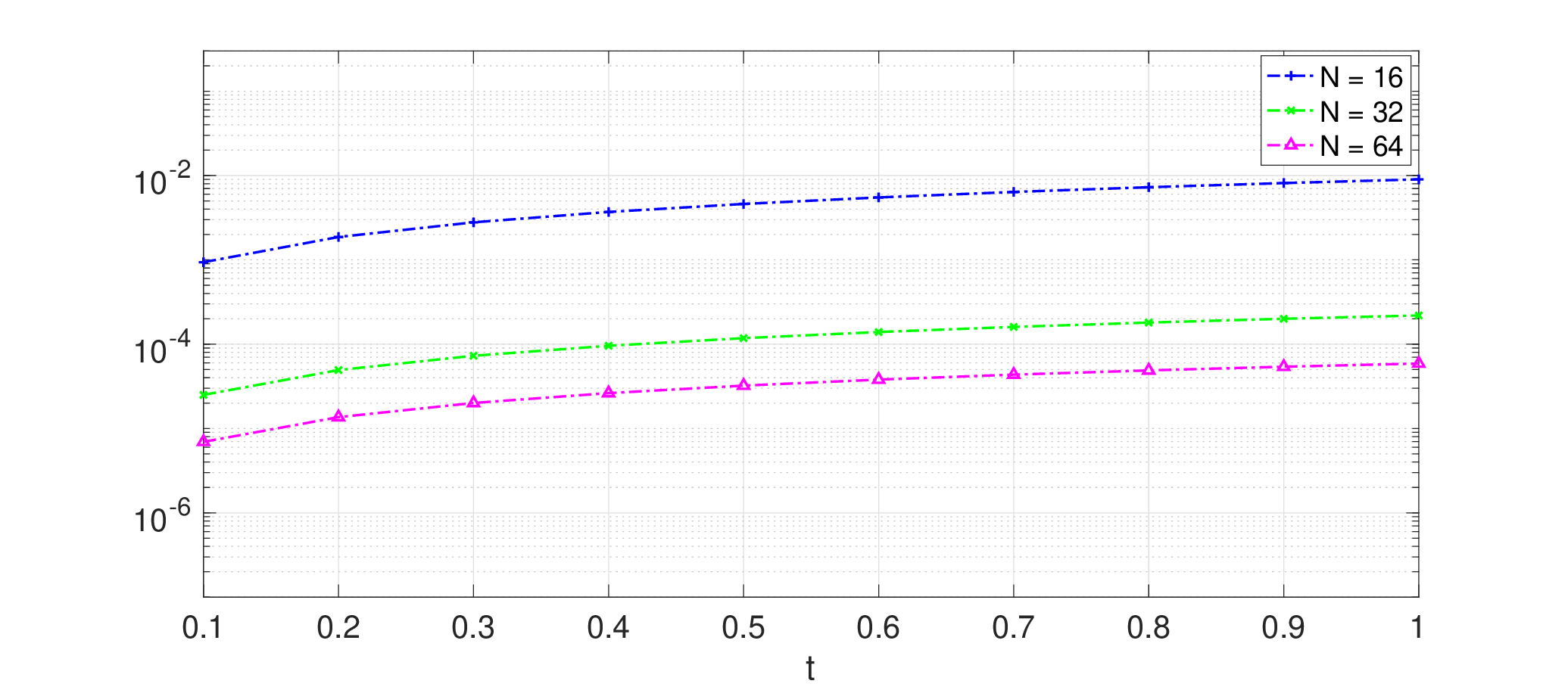}
            \caption{Energy conservation $|E(t) - E(0)|$ in 2D.}
            \label{fig-energy-2D}
        \end{subfigure}
        \caption{Evolution of mass and energy in 2D for $N = N_r, N_s = N_{sig} = 12$ and $\Delta t = 0.1$.}
        \label{fig:mass-energy-2D}
\end{figure}

\textit{Example 3} (Discontinuous case).
We consider the following discontinuous initial data in 2D:
\begin{equation}\label{disinitial}
f^0_{\text{dis}}(\bk)=
\begin{cases}
\frac{\rho_{1}}{2\pi T_1}\exp\left(-\frac{|\bk|^2}{2T_{1}}\right),\quad \text{for}\quad \bk_{1}>0,\\[6pt]
\frac{\rho_{2}}{2\pi T_2}\exp\left(-\frac{|\bk|^2}{2T_{2}}\right),\quad \text{for}\quad \bk_{1}<0,
\end{cases}
\end{equation}
where we pick $\rho_{1}=\frac{6}{5}$, $\rho_{2}=\frac{4}{5}$, $T_{1}=\frac{2}{3}$, $T_{2}=\frac{3}{2}$ such that
\begin{equation*}
\rho(0) = \int_{\mathbb{R}^{2}} f^0_{\text{dis}} \,\rd \bk =1,  \quad  E(0) = \int_{\mathbb{R}^{2}} f^0_{\text{dis}}|\bk|^2\,\rd \bk = 2.
\end{equation*} 
The profile of the discontinuous initial condition is drawn in Fig.~\ref{fig-disinitial}.
We use this example to check the conservation property of our method. This is shown in Fig.~\ref{fig-mass-2D} and Fig.~\ref{fig-energy-2D}, where we present the conservation of mass and energy over time. 

\begin{figure}[htp]
\begin{adjustwidth}{-0.5cm}{-0.5cm}
\centering
            \begin{subfigure}[t]{0.49\linewidth}
            \includegraphics[width=\textwidth]{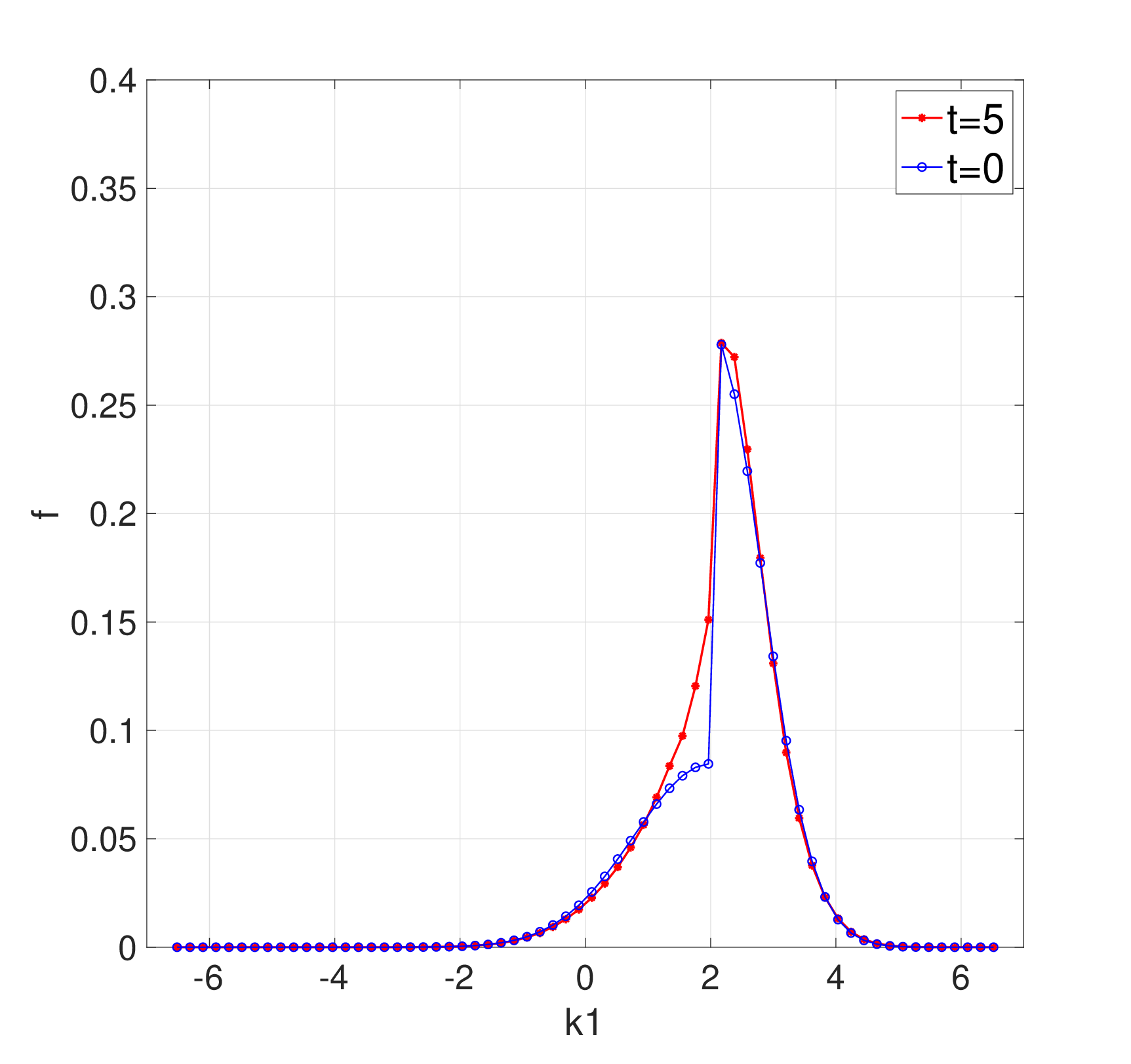}
            \caption{t=5}
        \end{subfigure}
        \begin{subfigure}[t]{0.49\linewidth}
            \includegraphics[width=\textwidth]{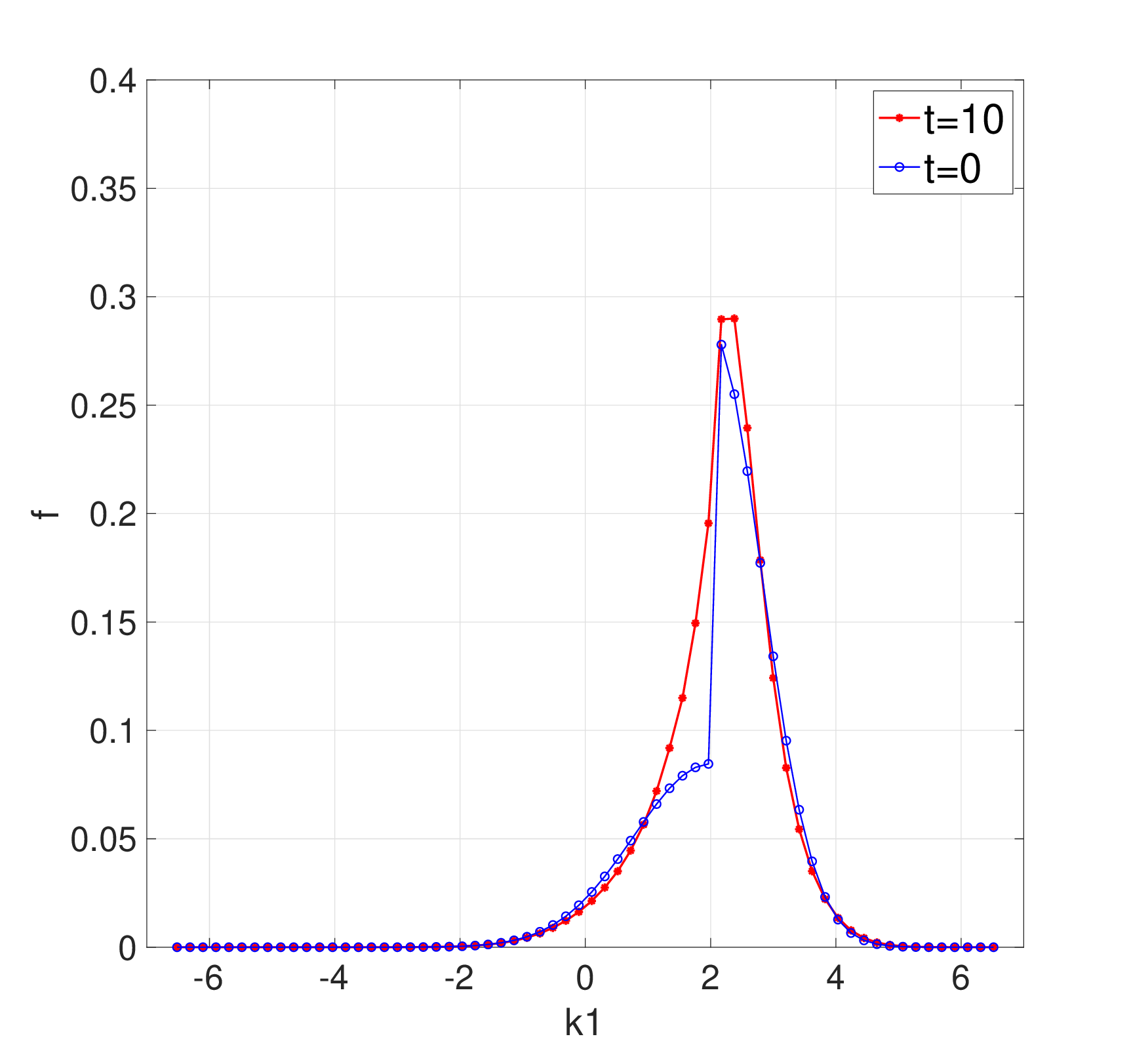}
            \caption{t=10}
        \end{subfigure}
        \begin{subfigure}[t]{0.49\linewidth}
            \includegraphics[width=\textwidth]{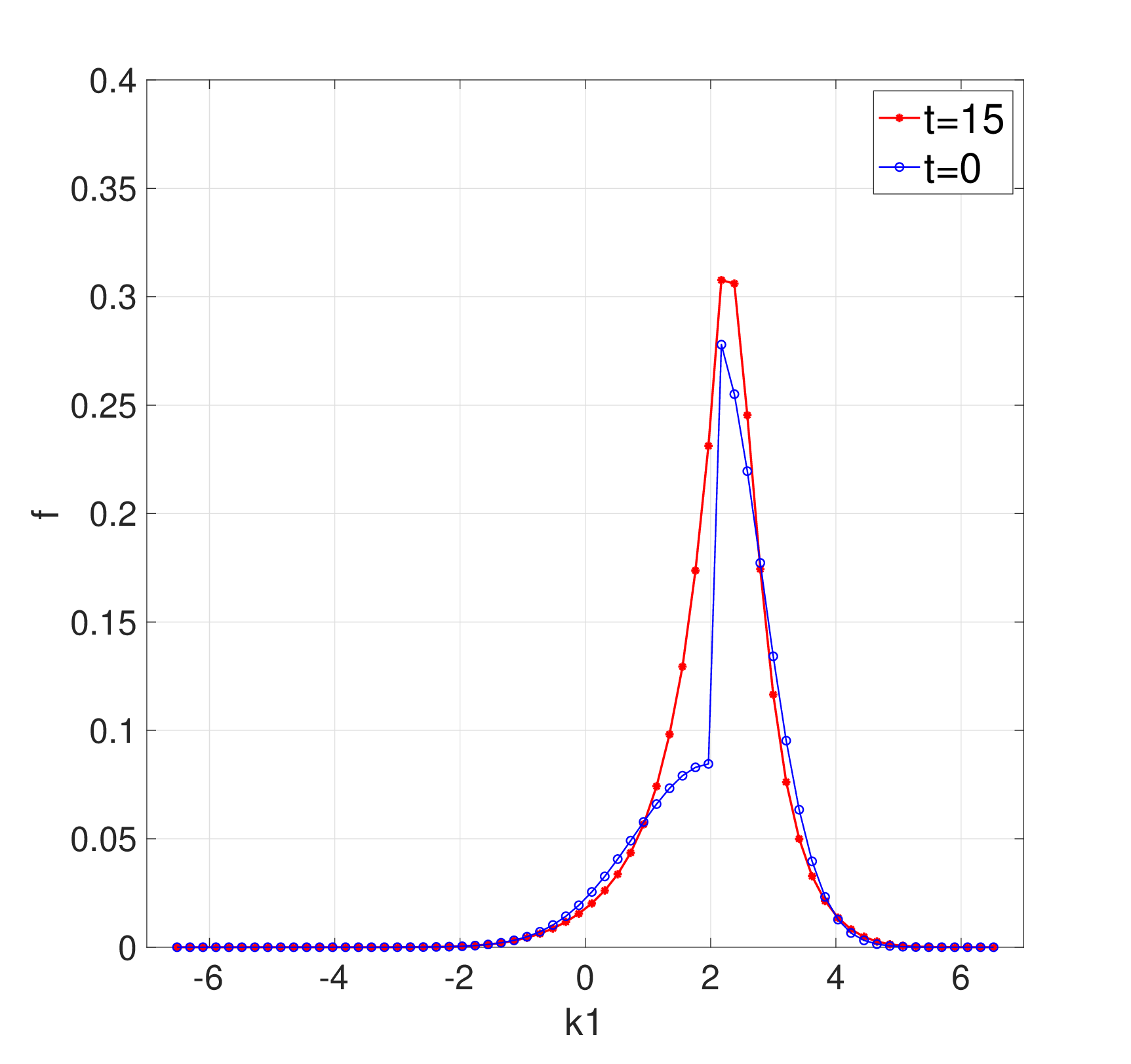}
            \caption{t=15}
        \end{subfigure}
        \begin{subfigure}[t]{0.49\linewidth}
            \includegraphics[width=\textwidth]{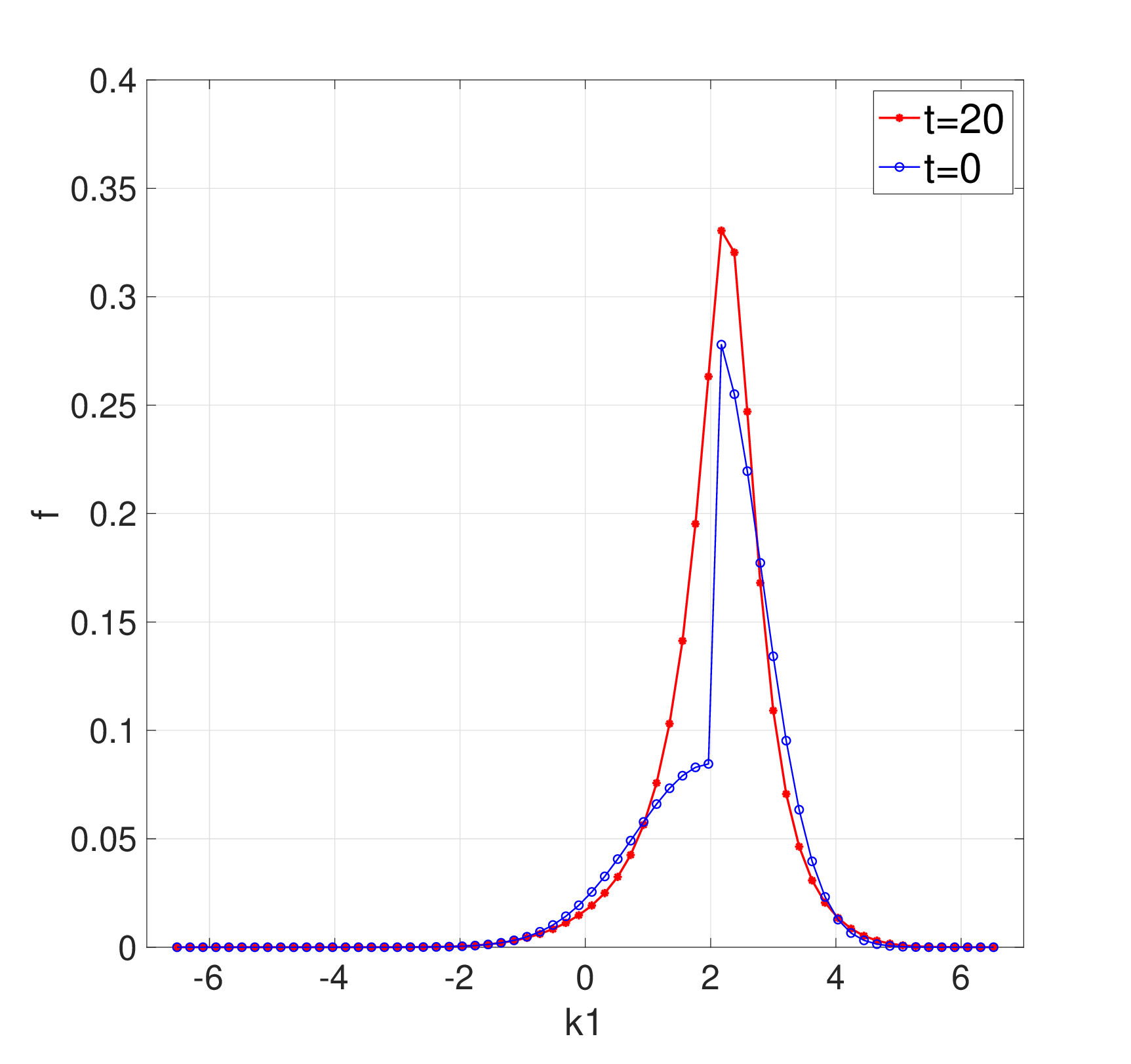}
            \caption{t=20}
        \end{subfigure}
        \newpage
        \begin{subfigure}[t]{0.49\linewidth}
            \includegraphics[width=\textwidth]{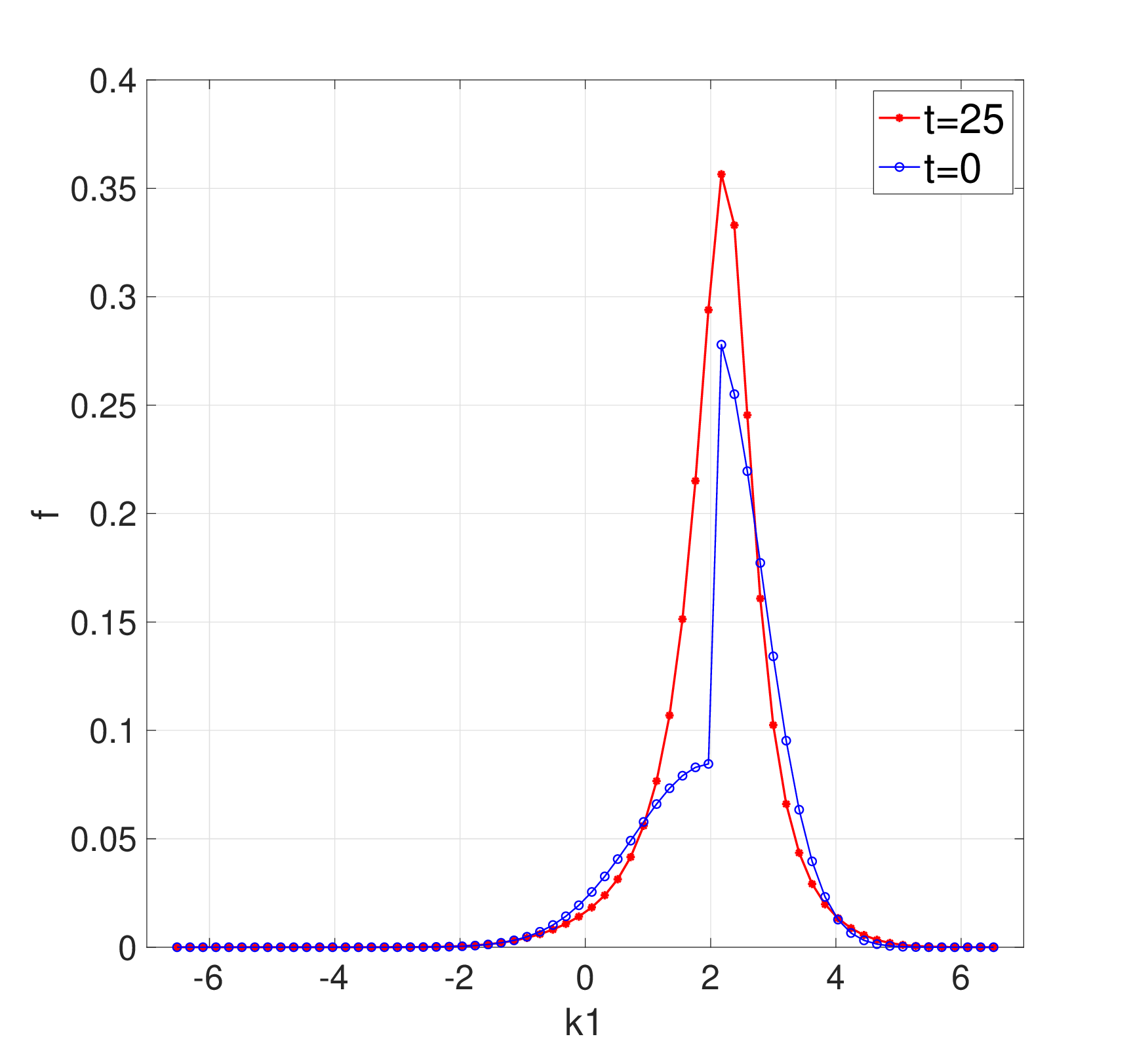}
            \caption{t=25}
        \end{subfigure}
        \begin{subfigure}[t]{0.49\linewidth}
            \includegraphics[width=\textwidth]{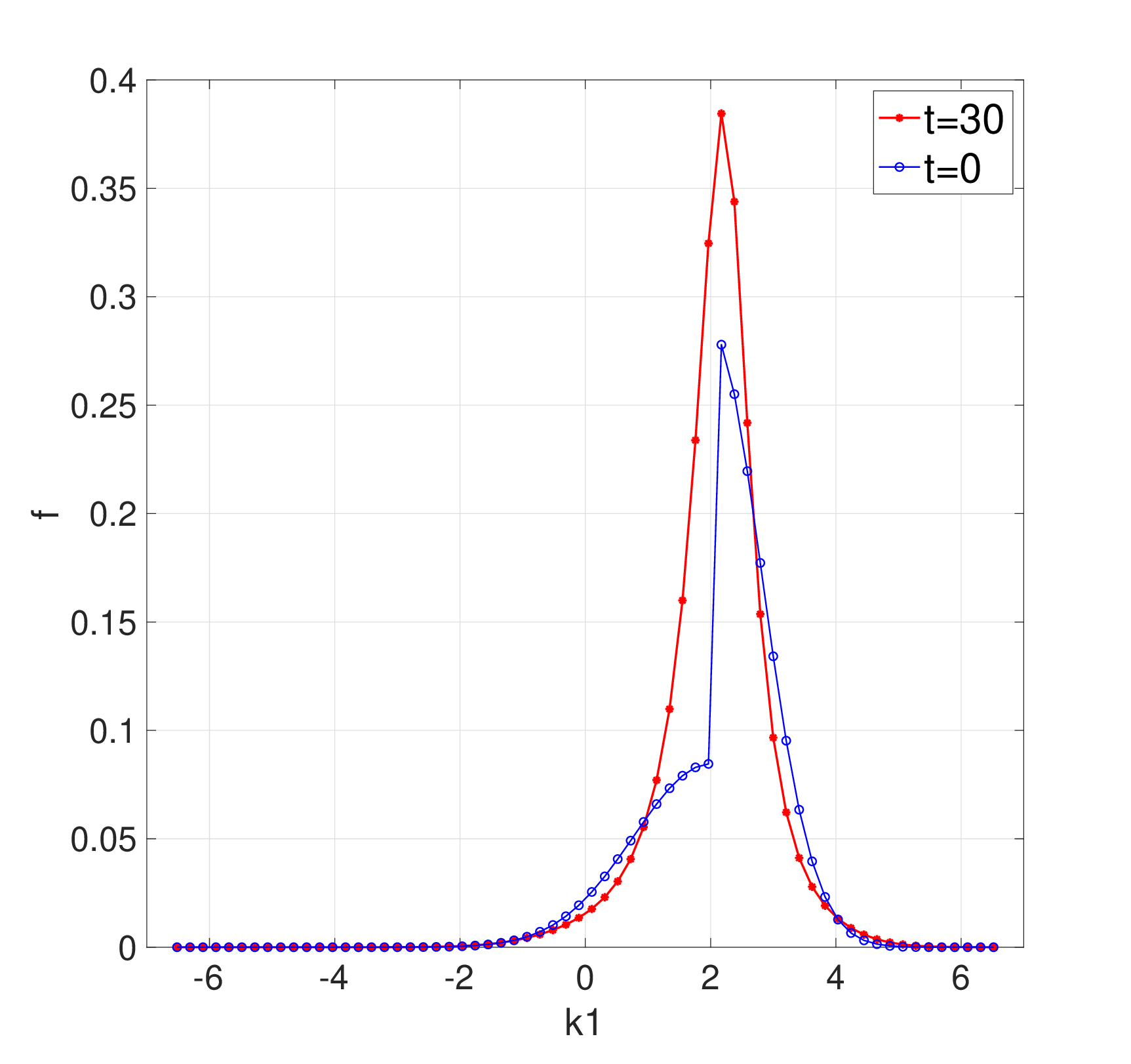}
            \caption{t=30}
        \end{subfigure}
\end{adjustwidth}
        \caption{Time evolution of the solution profile for the non-isotropic initial condition with $ N = N_r = 64, N_s = N_{sig} = 8$, $S=3$ and $\Delta t = 1$.}
        \label{fig-surf-aniso-dis}
    \end{figure}

\textit{Example 4} (Non-isotropic case). In this example, we test a conjecture regarding the non-isotropic case. As discussed in the theoretical work \cite[pp. 5-7]{EV2015-1} and demonstrated by the numerical experiment in Example 3, an isotropic initial condition can lead to a blow-up of the solution at the origin. It is conjectured that a non-isotropic initial condition will also lead to blow-up, but not symmetrically. In this test, we observe this phenomenon. As shown in Fig.~\ref{fig-surf-aniso-dis}, the solution evolves into a more concentrated profile; however, the center of the concentration does not remain at the origin, but shifts to a location dependent on the initial condition.

\subsection{3D case}
\label{subsec:3D}

\subsubsection{Stationary state test in 3D}
\label{subsubsec:stationary-3D}

Similar to the 2D case, we first use the stationary solution \eqref{equi} to verify the accuracy of the method in 3D.

\textit{Example 5} (Stationary solution in 3D). 
Consider the specific stationary solution \eqref{equi} in 3D with $\mu = 10$, $\bm{\nu} = (5,5)$ and $\xi=20$. The $L^{\infty}$ and $L^2$ errors  $\mathcal{K}(f_{\text{eq}})$  are presented in Table \ref{table-3D-Keq}, while the associated computation time for each components:  $\mathcal{K}_1$, $\mathcal{K}_2$, $\mathcal{K}_3$ and $\mathcal{K}_4$, are recorded in Table \ref{table-3D-time}.

\begin{table}[h!]
\begin{tabular}{| c || c | c |} 
 \hline
  $N$ & $\|\mathcal{K}(f_{\text{eq}})\|_{L^{\infty}}$ & $\| \mathcal{K}(f_{\text{eq}})\|_{L^{2}}$  \\ [0.5ex]
 \hline
 $16$ & $1.7066 \times 10^{-4}$ & $1.7286 \times 10^{-6}$ \\  
 \hline
 $32$ & $1.2351\times 10^{-4}$ & $1.9602\times 10^{-7}$ \\ 
 \hline
 $64$ & $7.2131\times 10^{-6}$ & $5.4791\times 10^{-8}$ \\ 
 \hline
\end{tabular}
\qquad
\begin{tabular}{| c || c | c |} 
 \hline
  $N$& $\|\mathcal{K}(f_{\text{eq}})\|_{L^{\infty}}$ & $\| \mathcal{K}(f_{\text{eq}})\|_{L^{2}}$  \\ [0.5ex]
 \hline
 $16$ & $1.9689 \times 10^{-4}$ & $7.2831 \times 10^{-7}$ \\  
 \hline
 $32$ & $6.8742\times 10^{-5}$ & $4.9058\times 10^{-8}$ \\ 
 \hline
 $64$ & $5.0208\times 10^{-6}$ & $9.9046\times 10^{-9}$ \\ 
 \hline
\end{tabular}
\caption{$L^{\infty}$ and $L^2$ error of stationary solution in 3D for $N_s = 6$ (left) and $N_s = 12$ (right) with $N_s = N_{sig}$ and $ S = 5$.}
\label{table-3D-Keq}
\end{table}
\begin{table}[h!]
\centering
\begin{tabular}{|c || c | c | c | c || c |} 
 \hline
   N & $\mathcal{K}_1$ & $\mathcal{K}_2$ & $\mathcal{K}_3$ & $\mathcal{K}_4$ &$\mathcal{K}$\\ [0.5ex]
 \hline
 $16$ & $0.27$s & $0.17$s & $0.20$s & $0.18$s & $0.82$s\\ 
 \hline
 $ 32$ & $2.11$s & $1.45$s & $1.60$s & $1.39$s & $6.55$s\\ 
 \hline
 $64$ & $22.38$s & $14.67$s & $17.73$s & $16.21$s & $70.99$s\\ 
 \hline
\end{tabular}
\caption{Computational time in 3D for $N_r = N, N_s = N_{sig} = 6$ and $S = 5$.}
\label{table-3D-time}
\end{table}

To further clarify the influence of $N_s$ and $N_{sig}$ in 3D, we compute $\mathcal{K}(f_{\text{eq}})$ and present the errors in $L^{\infty}$ and $L^{2}$ norms in Table.~\ref{table-3D-Keq-Nsig} for different $N_s$ and $N_{sig}$, which numerically demonstrate the validity of the approximation via the Spherical Design quadrature approach in the fast algorithm introduced in Sec.~\ref{subsec:fast}.

\renewcommand{\arraystretch}{1.2}
\begin{table}[h!]
\begin{tabular}{| c || c | c |} 
 \hline
  $N_s = N_{sig}$ & $\|\mathcal{K}(f_{\text{eq}})\|_{L^{\infty}}$ & $\| \mathcal{K}(f_{\text{eq}})\|_{L^{2}}$  \\ [0.5ex]
 \hline
 $6$ & $1.2351 \times 10^{-4}$ & $1.9602 \times 10^{-7}$ \\  
 \hline
 $12$ & $6.8742\times 10^{-5}$ & $4.9058\times 10^{-8}$ \\ 
 \hline
 $32$ & $3.4967\times 10^{-5}$ & $1.6784\times 10^{-8}$ \\ 
 \hline
\end{tabular}
\caption{$L^{\infty}$ and $L^2$ error of stationary solution in 3D for $N_s = N_{sig}$ with $N = N_{r} = 32$ and $ S = 5$.}
\label{table-3D-Keq-Nsig}
\end{table}


\subsubsection{Time evolution test in 3D}
\label{subsubsec:time-3D}

We then solve the time dependent problem \eqref{WKE} in 3D. Time discretization is again performed using the fourth-order Runge-Kutta (RK4) method.

\begin{figure}[h!]
        \centering
        \begin{subfigure}[t]{1.1\textwidth}
            \includegraphics[width=1\textwidth]{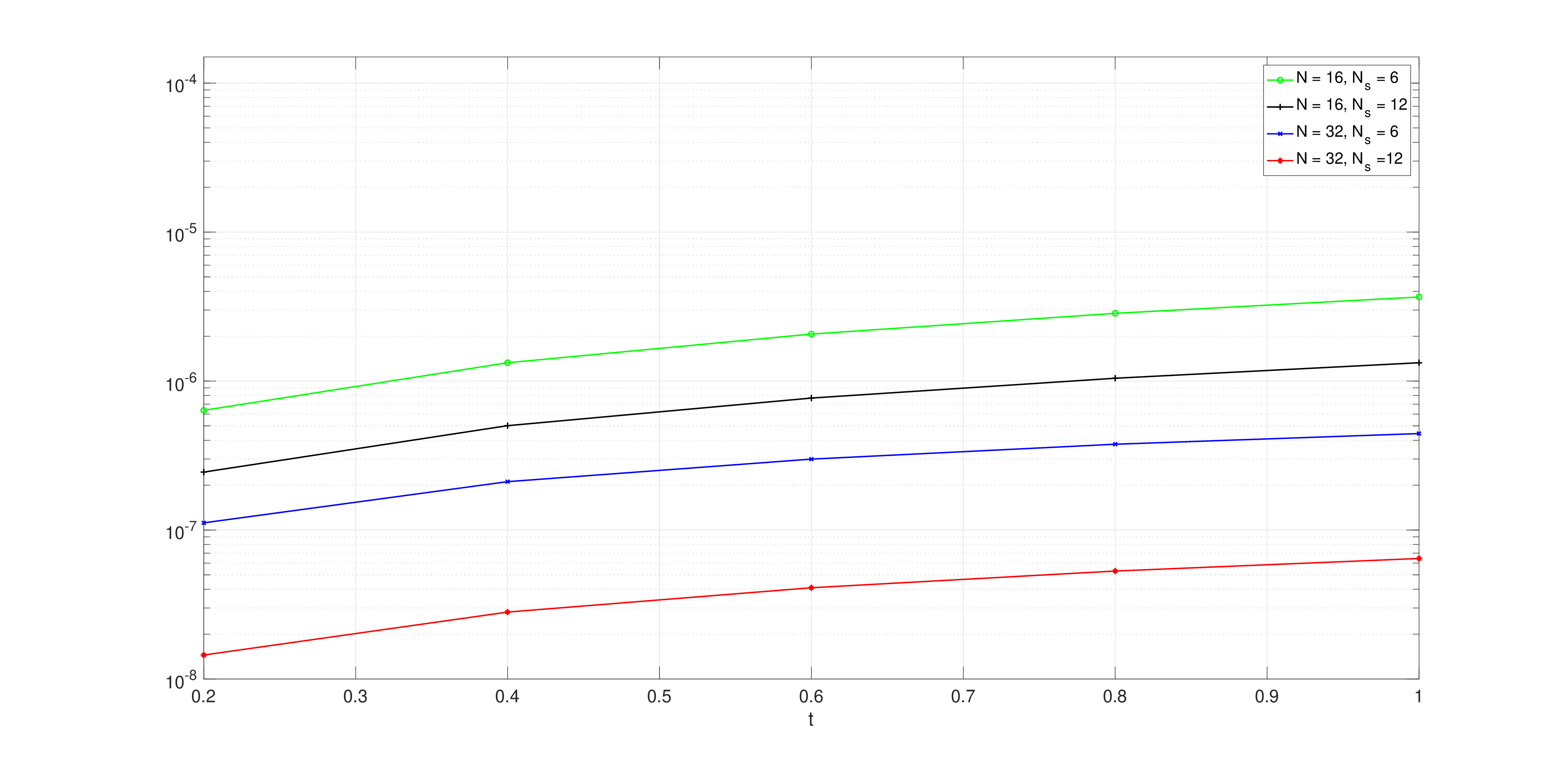}
            \caption{Mass conservation $|\rho(t) - \rho(0)|$ in 3D}
            \label{fig-mass-3D}
        \end{subfigure}
        \begin{subfigure}[t]{1.1\textwidth}
            \includegraphics[width=1\textwidth]{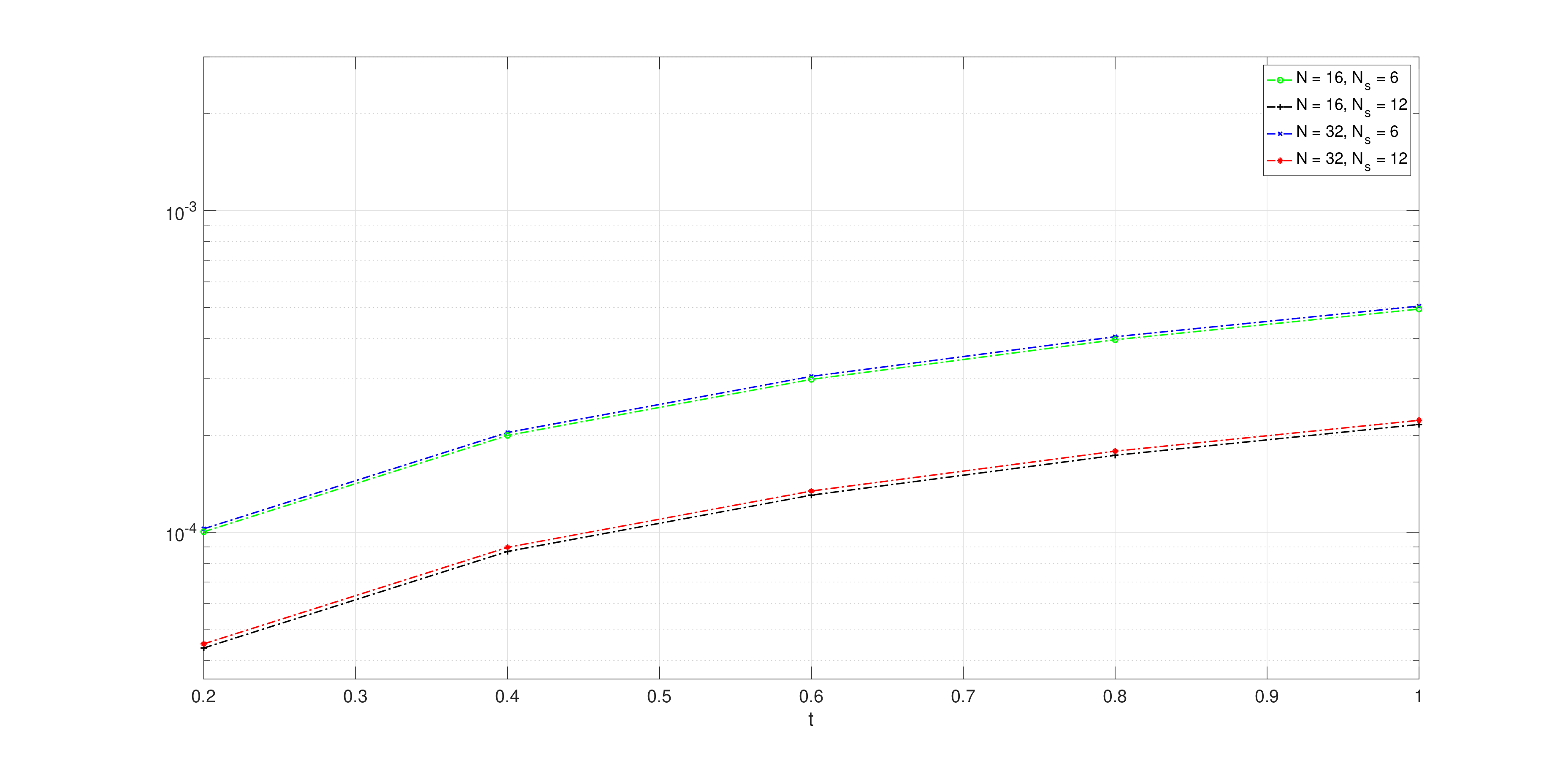}
            \caption{Energy conservation $|E(t) - E(0)|$ in 3D}
            \label{fig-energy-3D}
        \end{subfigure}
        \caption{Evolution of mass and energy in 3D for $ N = N_r, N_s = N_{sig}$, $S=0.33$ and $\Delta t= 0.1$.}
        \label{mass-energy-3D}
\end{figure}

\textit{Example 6} (Mass and energy conservation in 3D).
Here we consider the same discontinuous initial condition \eqref{disinitial} in 3D. The mass and energy evolution is presented in Fig.~\ref{mass-energy-3D}.

\section{Conclusion}
\label{sec:conclusion}
In this work, we develop a fast Fourier spectral method for solving the wave kinetic equation (WKE), which arises as the kinetic limit of the cubic Schr{\"o}dinger equation. The core idea of our approach is to reformulate the high-dimensional nonlinear wave kinetic operator as a spherical integral, drawing an analogy to the classical Boltzmann collision operator. Additionally, we extract a double-convolution structure from the numerical system by applying the Fourier spectral approximation to the solution, enabling the use of the fast Fourier transform (FFT) to significantly enhance computational efficiency. This FFT-based approach not only accelerates computation but also eliminates the need for pre-computation and additional storage, making it a practical and scalable solution for high-dimensional problems.

While our approach borrows key ideas from previous works \cite{MP06, HY15}, its application to the wave kinetic equation (WKE) is novel and could have a significant impact. Compared to earlier methods for computing the WKE, such as the Discrete Interaction Approximation (DIA) \cite{hasselmann1985computations} and quadrature-based approaches (e.g., the WRT method) \cite{webb1978non, tracy1982theory}, our method excels in both accuracy and efficiency.
Through a series of numerical experiments in 2D and 3D, we demonstrate the performance of the proposed method in both steady-state and time-evolving scenarios, including isotropic and non-isotropic cases. These results highlight the potential of our method as a robust and efficient tool for simulating wave kinetic dynamics within wave turbulence theory.

\section*{Acknowledgments}
We would like to thank Yu Deng and Angeliki Menegaki for their helpful discussion about the wave kinetic equation. K.~Qi also would like to thank Jingwei Hu for her helpful guidance and discussion about the spectral method for Boltzmann type equation, and thank Zaher Hani for his helpful discussion about wave kinetic equation and introduction of reference \cite{BBKKS2022} during the 2024 Rivi\`ere-Fabes Symposium at the University of Minnesota.

\appendix
\section{Reformulation of WKE in the form of classical Boltzmann-type equation}
\label{appen:derivation}

To make the paper self-explained, in this appendix, we recall the derivation that reformulates the form \eqref{K} as an integral over a sphere, leading to the usual form of Boltzmann-type collision operator as in \eqref{WKE-Boltz}.

Step I: Denote $ \mathbf{O} = \frac{(\bk_1 + \bk_3)}{2}$ and $\mathbf{P} = \frac{(\bk_1 - \bk_3)}{2}$ such that
\begin{equation*}
    \bk_1 = \mathbf{O} + \mathbf{P}, \quad \bk_3 = \mathbf{O} - \mathbf{P},
\end{equation*}
and applying the change of variable $(\bk_1, \bk_3)\mapsto (\mathbf{O}, \mathbf{P})$ with the Jacobian $2^{d}$ in \eqref{K}, we have
\begin{equation*}
    \begin{split}
      \mathcal{K}(f,f,f)(\bk) = 2^{d}&\int_{(\RR^d)^3}  \Big[ f(\mathbf{O}+\mathbf{P}) f(\bk_2) f(\mathbf{O}-\mathbf{P}) - f(\bk) f(\bk_2) f(\mathbf{O}-\mathbf{P})\\
      &\quad+ f(\mathbf{O}+\mathbf{P}) f(\bk) f(\mathbf{O}-\mathbf{P})
      - f(\mathbf{O}+\mathbf{P}) f(\bk_2) f(\bk)\Big]\\
    &\times \delta\left(\bk_2 + \bk - 2\mathbf{O} \right) \times \delta \left( \frac{|\bk_2|^2}{2} + \frac{|\bk|^2}{2} - \mathbf{O}^2 - \mathbf{P}^2 \right) \,\rd \mathbf{O} \,\rd \bk_2 \,\rd \mathbf{P}.
    \end{split}
\end{equation*}

Step II: Let $\mathbf{P} = \rho \bsigma$ with $\rho \in [0,\infty)$ and $\bsigma \in \mathbb{S}^{d-1}$, such that $\rd \mathbf{P} = \rho^{d-1} \,\rd \rho \,\rd \bsigma$,
\begin{equation*}
    \begin{split}
      &\mathcal{K}(f,f,f)(\bk) =\\[4pt]
      & 2^{d}\int_{(\RR^d)^2} \int_{\mathbb{S}^{d-1}} \int_0^{\infty}  \Big[ f(\mathbf{O} +\rho \bsigma ) f(\bk_2) f(\mathbf{O} -\rho \bsigma) - f(\bk) f(\bk_2) f(\mathbf{O}-\rho \bsigma) \\
      &\quad + f(\mathbf{O}+\rho \bsigma) f(\bk) f(\mathbf{O}-\rho \bsigma) - f(\mathbf{O}+\rho \bsigma) f(\bk_2) f(\bk)\Big]\\
    &\times \delta\left(\bk_2 + \bk - 2\mathbf{O} \right) \times \delta \left( \frac{|\bk_2|^2}{2} + \frac{|\bk|^2}{2} - \mathbf{O}^2 - \rho^2 \right) \rho^{d-1} \,\rd \rho \,\rd \bsigma \,\rd \mathbf{O} \,\rd \bk_2,
    \end{split}
    \end{equation*}
    then, replacing $\rho^2$ by $= \frac{|\bk_2|^2}{2}+\frac{|\bk|^2}{2}-\mathbf{O}^2$, we obtain 
    \begin{equation*}
    \begin{split}
    \mathcal{K}&(f,f,f)(\bk) =
    \frac{1}{2}\int_{(\RR^d)^2} \int_{\mathbb{S}^{d-1}} \\ 
    &\Bigg[ f\left(\mathbf{O}+\sqrt{\frac{|\bk_2|^2}{2}+\frac{|\bk|^2}{2}-\mathbf{O}^2} \bsigma\right) f(\bk_2) f\left(\mathbf{O}-\sqrt{\frac{|\bk_2|^2}{2}+\frac{|\bk|^2}{2}-\mathbf{O}^2} \bsigma\right) \\
    &- f(\bk) f(\bk_2) f\left(\mathbf{O}-\sqrt{\frac{|\bk_2|^2}{2}+\frac{|\bk|^2}{2}-\mathbf{O}^2} \bsigma\right) \\
    &+ f\left(\mathbf{O}+\sqrt{\frac{|\bk_2|^2}{2}+\frac{|\bk|^2}{2}-\mathbf{O}^2} \bsigma \right) f(\bk) f\left(\mathbf{O}-\sqrt{\frac{|\bk_2|^2}{2}+\frac{|\bk|^2}{2}-\mathbf{O}^2} \bsigma\right) \\
    &- f\left(\mathbf{O}+\sqrt{\frac{|\bk_2|^2}{2}+\frac{|\bk|^2}{2}-\mathbf{O}^2} \bsigma \right) f(\bk_2) f(\bk)\Bigg]\\
    &\times \delta\left(\mathbf{O} - \frac{\bk_2 + \bk}{2}  \right) \times   \left(\frac{|\bk_2|^2}{2}+\frac{|\bk|^2}{2} - \mathbf{O}^2\right)^{\frac{d-2}{2}} \,\rd \bsigma \,\rd \mathbf{O} \,\rd \bk_2.
    \end{split}
    \end{equation*}
    
Step III: By integration over $\mathbf{O}$ (replacing $\mathbf{O}$ by $= \frac{\bk_2 + \bk}{2}$ above) and considering the fact that
    \begin{equation*}
        \sqrt{\frac{|\bk_2|^2}{2}+\frac{|\bk|^2}{2}-\mathbf{O}^2}  = \frac{|\bk-\bk_2|}{2},
    \end{equation*}
we obtain
    \begin{equation*}
    \begin{split}
    \mathcal{K}(f,f,f) =& \frac{1}{2}\int_{\RR^d} \int_{\mathbb{S}^{d-1}}  \Bigg[ f\left(\frac{\bk_2 + \bk}{2} + \frac{|\bk-\bk_2|}{2}\bsigma\right) f(\bk_2) f\left(\frac{\bk_2 + \bk}{2} - \frac{|\bk-\bk_2|}{2}\bsigma\right) \\[3pt]
    &- f(\bk) f(\bk_2) f\left(\frac{\bk_2 + \bk}{2} - \frac{|\bk-\bk_2|}{2}\bsigma\right) \\[3pt]
    &+ f\left(\frac{\bk_2 + \bk}{2} + \frac{|\bk-\bk_2|}{2}\bsigma\right) f(\bk) f\left(\frac{\bk_2 + \bk}{2} - \frac{|\bk-\bk_2|}{2}\bsigma\right) \\[3pt]
    &- f\left(\frac{\bk_2 + \bk}{2} + \frac{|\bk-\bk_2|}{2}\bsigma\right) f(\bk_2) f(\bk) \Bigg]\\[3pt]
    &\times \left(\frac{|\bk-\bk_2|}{2}\right)^{d-2} \,\rd \bsigma \,\rd \bk_2,
    \end{split}
\end{equation*}
which, by denoting
\begin{equation*}
   \frac{\bk_2 + \bk}{2} + \frac{|\bk-\bk_2|}{2}\bsigma =:\bk', \quad 
    \frac{\bk_2 + \bk}{2} - \frac{|\bk-\bk_2|}{2}\bsigma = : \bk'_2 \,,
\end{equation*}
can finally be written as in \eqref{WKE-Boltz}.

\section{Mass conservation of the spectral numerical system}
\label{appen:mass}
In this appendix, we show that the numerical system without using the fast algorithm preserves mass, that is,
\begin{equation} 
\int_{\D} f_N(t,\bk) \,\rd \bk=\int_{\D} f^{0}_N(\bk) \,\rd \bk.
\end{equation}
where $f^{0}_N(\bk)$ is the approximation of the initial condition $f^{0}(\bk)$ as in \eqref{fj}.
Note that 
\begin{equation} 
\int_{\mathcal{D}_L} f_N(t,\bk)\,\rd{\bk} = \sum_{|\bj|=-\frac{N}{2}}^{\frac{N}{2}-1} \hat{f}_{\bj}(t) \int_{\mathcal{D}_L}  \e^{\im \frac{\pi}{L}\bj\cdot \bk} \,\rd{\bk}=(2L)^d \hat{f}_{\bj=\mathbf{0}}(t),
\end{equation}
where $\hat{f}_{\bj=\mathbf{0}}(t)$ is the zero-th mode of the numerical solution and is governed by
\begin{equation} 
\label{eq:fbj0}
\partial_{t} \hat{f}_{\bj = \mathbf{0}} = \hat{\mathcal{K}}_{\bj = \mathbf{0}}^R.
\end{equation}
By applying the weak form \eqref{weak-WKE} to \eqref{WKE-Boltz} and choosing the test function $\varphi(\bk) = \e^{-\im \frac{\pi}{L} \bj \cdot \bk}$, the numerical evaluation of $\hat{\mathcal{K}}_{\bj}^R$ in \eqref{Kj} becomes 
\begin{equation}
    \begin{split}
        &\hat{\mathcal{K}}_{\bj}^R
        = \frac{1}{(2L)^{d}} \sum_{|\bl|, |\m|, |\n| = -\frac{N}{2}}^{\frac{N}{2}-1} \int_{\D_L} \int_{\B_R} \int_{\mathbb{S}^{d-1}} |\bq|^{d-2} \hat{f}_{\bl} \hat{f}_{\m} \hat{f}_{\n}\\
        & \times \Big( \e^{\im \frac{\pi}{L} (\bl+\m+\n-\bj) \cdot \bk}
          \e^{-\im \frac{\pi}{L} \frac{1}{2} (2\m+\bl+\n) \cdot \bq} 
          \e^{\im \frac{\pi}{L} \frac{1}{2}|\bq|(\bl-\n)\cdot\bsigma} - \e^{\im \frac{\pi}{L} \bl \cdot \bk}
          \e^{\im \frac{\pi}{L} \m \cdot (\bk-\bq)} \e^{\im \frac{\pi}{L} \n \cdot [\bk - \frac{1}{2}(\bq + |\bq|\bsigma)]} \\
          & + \e^{\im \frac{\pi}{L} (\bl+\m+\n-\bj) \cdot \bk}
          \e^{-\im \frac{\pi}{L} \frac{1}{2} (\bl+\n) \cdot \bq} 
          \e^{\im \frac{\pi}{L} \frac{1}{2}|\bq|(\bl-\n)\cdot\bsigma} - \e^{\im \frac{\pi}{L} (\bl+\m+\n-\bj) \cdot \bk}
          \e^{-\im \frac{\pi}{L} \frac{1}{2} (2\m+\bl) \cdot \bq} 
          \e^{\im \frac{\pi}{L} \frac{1}{2}|\bq|\bl\cdot\bsigma}  \Big)\\
        &\times \Big[ \e^{-\im \frac{\pi}{L} \bj \cdot \bk} + \e^{-\im \frac{\pi}{L} \bj \cdot \bk_2} - \e^{-\im \frac{\pi}{L} \bj \cdot \bk'} - \e^{-\im \frac{\pi}{L} \bj \cdot \bk'_2} \Big]
          \,\rd \bsigma \,\rd \bq  \,\rd \bk
    \end{split}
\end{equation}
and this form can be easily verified that $\hat{\mathcal{K}}_{\bj }^R \equiv0$ when $\bj=\mathbf{0}$. Hence, it further implies from \eqref{eq:fbj0} that $\hat{f}_{\bj = \mathbf{0}}$ remains constant as time evolves. However, it is worth mentioning that the mass might not be strictly conserved if applying the fast algorithm, where the weight functions are approximated as in \eqref{315} and the integrals in $\bk, \sigma, \bq$ are replaced by a quadrature sum.

\bibliographystyle{amsplain}
\bibliography{Qi_bib}

\end{document}